\documentclass[pdf14]{statsoc}

\usepackage[a4paper,margin=1in,bottom=1in,top=0.9in]{geometry}

\usepackage{amsmath}
\usepackage{amsfonts}
\usepackage{amssymb}
\usepackage{graphicx}
\usepackage{color, graphics}
\usepackage{natbib} 
\usepackage{bm}

\usepackage{subcaption}

\newcommand{\iidsim}{\overset{i.i.d.}{\sim}}
\newcommand{\indsim}{\overset{indep}{\sim}}
\newcommand{\eqdist}{\overset{d}{=}}

\newcommand{\approxiid}{\overset{iid}\approx}

\newtheorem{thm}{Theorem}[section]
\newtheorem{lemma}%[thm]
{Lemma}[section]
\newtheorem{prop} %[thm]
{Proposition}[section]
{Corollary}[section]
\newtheorem{defn} %[thm]
{Definition}[section]
\newtheorem{rmrk} %[thm]
{Remark}[section]
{Example}[section]

\newcommand{\Polya}{P\'{o}lya }
\newcommand{\Norm}{\text{N}} %Normal distribution
\newcommand{\DP}{\text{DP} } %Dirichlet Process

\newcommand{\bpi}{\mathbf{\pi}}

\newcommand{\tG}{\tilde{G}}
\newcommand{\tP}{{\tilde{P}}}
\newcommand{\tF}{{\tilde{F}}}
\newcommand{\cF}{{\cal F}}

 \title[Quasi-Bayes properties]{Quasi-Bayes properties of a procedure for sequential learning in mixture models}
 \author{Sandra Fortini}
 \address{Bocconi University, Milan, Italy.}
 \author[S. Fortini, S. Petrone]{Sonia Petrone}
 \address{Bocconi University, Milan, Italy.}
 \email{sonia.petrone@unibocconi.it}
 
\begin{document}
 \maketitle

%% tagliare questo sotto se NON JRSS
%\title{Quasi-Bayes properties of a procedure for sequential learning in mixture models}
%\author{Sandra Fortini and Sonia Petrone}
%\begin{document}
%\maketitle
%%%

\begin{abstract}

Bayesian methods are often optimal, yet increasing pressure for fast computations, especially with streaming data, brings renewed interest in faster, possibly sub-optimal, solutions. The extent to which these algorithms approximate Bayesian solutions is a question of interest, but often unanswered.  
We  propose a methodology to address this question in predictive settings, when the algorithm can be reinterpreted as a probabilistic predictive rule. 
We specifically develop the proposed methodology for a recursive procedure
for online learning in nonparametric mixture models, often refereed to as Newton's algorithm. This algorithm is simple and fast; however, its approximation properties are unclear. 
By reinterpreting it as a predictive rule, we can show that it underlies a  statistical model which is,  asymptotically, a Bayesian, exchangeable mixture model. In this sense, the recursive rule provides a quasi-Bayes solution. 
While the algorithm only offers a point estimate, our clean statistical formulation
allows us to provide the  asymptotic posterior distribution and asymptotic credible intervals for the mixing distribution. Moreover, 
it gives insights for tuning the parameters, as we illustrate in simulation studies, and paves the way to  extensions in various directions. Beyond mixture models, our approach can be applied to other predictive algorithms. 

\end{abstract}


\section{Introduction}

Bayesian methods have always been attractive, for their internal coherence, their rigorous way of quantifying uncertainty through probability and their optimal properties in many problems. Analytic difficulties have been overcome by  efficient computational methods and Bayesian procedures are nowadays widely and successfully used in many fields. However, fast computations remain a challenge, that hampers an even wider application of Bayesian methods among practitioners, the more so with streaming data and online learning, where inference and prediction have to be continuously updated as new data become available. 
In the modern trade off between statistical and computational efficiency, slightly misspecified but computationally more tractable methods receive renewed interest, as a reasonable compromise. Popular algorithms, such as Approximate Bayesian Computation (ABC) and variational Bayes, arise as approximations of an optimal Bayesian solution. Indeed, one could expect 
 that a method which performs well is at least {\em approximately Bayes}. For a Bayesian statistician, the capacity of a learning scheme to be, at least approximately, a Bayesian learning scheme should be a minimal requirement for its validation.

We  propose a methodology to address the above questions in predictive settings, when the  
algorithm can be reinterpreted as a probabilistic predictive rule, that implicitly defines an underlying statistical model. We then leverage on characterizing properties of the predictive rule to obtain such model explicitly.  This approach allows to  develop the algorithm into a clean statistical procedure and, on this basis, clarify its properties as an approximation of a fully Bayesian  method. 
Predictive constructions are a powerful tool in Bayesian inference, to characterize prior laws;  %let us just mention  \cite{blackwellMacqueen1973} \Polya sequences(or Chinese Restaurant Process), that give a predictive characterization of the Dirichlet process. Using predictive characterizations 
yet, their use in the problem under study appears novel.

% moreover, we aim at at finding the statistical model, besides the prior law. 
%that we will say to be a quasi-Bayes approximation of a given Bayesian procedure if the underlying model is equivalent to the Bayesian model, at least asymptotically. 

We specifically develop the proposed predictive methodology in the important case of sequential learning in mixture models. There is an extensive literature on Bayesian learning for mixtures. However, sequential learning, specifically on the mixing distribution, on which we focus, is less developed. Moreover, most popular Bayesian nonparametric mixture models, e.g. Dirichlet process mixtures, assume a discrete mixing distribution. The case of an absolutely continuous mixing distribution, important, for example in multiple shrinkage estimation \citep{george1986}, is also less developed. \cite{petroneVeronese2002} use a general extension of the Bernstein polynomial prior on the latent distribution, but computations require MCMC and estimation with streaming data is not addressed. 

An interesting recursive procedure for unsupervised sequential learning and classification in finite mixtures was proposed by \cite{smithMakov1978} and extended by M. Newton and collaborators (\cite{newtonQuintanaZhang1998}, \cite{newtonZhang1999},  \cite{quintanaNewton2000}, \cite{newton2002}) to provide a fast, approximately Bayes, solution  in nonparametric mixture models. A thoughtful review is given by \cite{martin2019}. 
Recent interesting developments are provided in \cite{hahnMartinWalker2018}. Convergence results have validated the recursive algorithm as a consistent frequentist estimator.
%(\cite{newtonZhang1999}, \cite{martinGhosh2008}, \cite{ghoshTokdar2006}, \cite{tokdarMartinGhosh2009}), 
%and recent work includes 
Further properties are given in 
\cite{favaroGuglielmiWalker2012} and \cite{zuanettiMueller2019}. However, the extent to which the recursive algorithm  provides an approximation of a Bayesian procedure is not fully understood. 

We  aim at shedding light on this question, by clarifying, through a predictive approach,  the statistical model underlying the recursive rule.  
This makes users aware of the assumptions implicitly made on the data when using the  algorithm, and of the related uncertainty. The proposed approach may be of interest as a method for quantifying the uncertainty of other predictive algorithms, beyond mixture models. 

Let us start with a first example, sequential unsupervised learning and classification by mixtures, considered by \cite{smithMakov1978}. 
The aim is to recursively classify observations $x_1, x_2, \ldots$ in one of $k$ populations (e.g., pattern types, or  signal sources, etc.), with no feedback about correctness of previous classifications. A finite mixture model for this task assumes 
\begin{equation}
\label{eq:finiteMixture}X_n \mid \bpi \iidsim \sum_{j=1}^k \pi_j \; f_j(x), 
\end{equation}
where i.i.d. stands for independent and identically distributed. Here the mixture components $f_j(\cdot)$ are known (extensive studies may be available on the specific components), but the mixing proportions $\bpi=(\pi_1, \ldots, \pi_k)$ are unknown. The classical Bayesian solution assigns a Dirichlet prior distribution the unknown proportions $\pi$, 
and proceeds by Bayes rule. Learning is solved through the posterior distribution $p(\bpi \mid x_1, \ldots, x_n)$ and classification through the predictive probabilities that $X_{n+1} \sim f_j$, given $(x_1, \ldots, x_n)$, for $j=1, \ldots, k$. Unfortunately, sequential computations are involved.

%Estimation of a latent density %and multiple shrinkage estimation, with streaming data is a related problem.  
The finite mixture model (\ref{eq:finiteMixture}) is a special case, 
for a discrete $\tG$ with atoms $1,\ldots, k$ having unknown masses $\pi_1, \ldots, \pi_k$, of a general mixture model 
\begin{equation} \label{eq:mixture}
 X_i \mid \tG \iidsim f_{\tG}(x)=\int f(x \mid \theta) d \tG(\theta).    
 \end{equation}
A problem of interest is to recursively estimate the latent distribution $\tG$ as new observations become available. 
In a Bayesian nonparametric approach, a prior with large support is assigned to the random mixing distribution $\tG$, a popular choice being a Dirichlet Process (DP), with parameters $\alpha$ and $G_0$, $\tG \sim \DP(\alpha, G_0)$. Then one proceeds by computing the conditional distributions of interest. Computations are involved, but can be addressed by MCMC methods, or via variational Bayes (\cite{blei2017}) or ABC approximations. 
If the observations $x_i$ arrive sequentially, one may resort to sequential Monte Carlo methods, sequential importance sampling (\cite{macEachernClydeLiu1999}), or more recent sequential variational Bayes methods (\cite{lin2013}, \cite{broderick2013}), or combinations of them (\cite{naesseth2018}). 
Still, these methods have a computational cost (for example, in the optimization steps) or can be derived only heuristically. The search for simple and fast recursive algorithms remains attractive. 
%Moreover, the discrete nature of the DP, and of most of its extensions, while extremely powerful in parsimoniously accounting for random partitions and clustering, is not appropriate for shrinkage problems, where the latent distribution should be absolutely continuous.  

The recursive rule proposed by M. Newton and collaborators, often referred to as Newton's algorithm,  starts at an initial guess $G_0$ and, for any $A$, recursively computes the estimated mixing distribution as 
\begin{equation} \label{eq:newton} 
G_n(A) = (1-\alpha_n) G_{n-1}(A) + \alpha_n 
\frac{\int_A f(x_n \mid \theta_n) d G_{n-1}(\theta_n)}
{ \int_\Theta f(x_n \mid \theta_n) dG_{n-1}(\theta_n)},   
\end{equation}
where $(\alpha_n)$ is a sequence of real numbers in $(0,1)$ and it is usually  assumed that $\alpha_n\rightarrow 0$ as $n\rightarrow\infty$, with $\sum_n\alpha_n=\infty$ and $\sum_n\alpha_n^2<\infty$. A standard choice, in analogy with DP mixtures, is $\alpha_n=1/(\alpha+n)$ with $\alpha>0$.  
For finite mixtures, as in (\ref{eq:finiteMixture}),  the rule (\ref{eq:newton}) corresponds to the sequential procedure of Smith and Makov (1978). 
 If $G_0$ has density $g_0$ with respect to a measure $\lambda$, then $G_n$ has density $g_n$ with respect to the same $\lambda$ and (\ref{eq:newton}) implies that $g_n$ satisfies the recursive equation
\begin{equation}\label{eq:recursiveDensities}
\begin{aligned}
g_n(\theta)&=(1-\alpha_n)g_{n-1}(\theta)+\alpha_n\frac{g_{n-1}(\theta)f(x_n\mid\theta)}{\int f(x_n\mid\theta')g_{n-1}(\theta')d\lambda(\theta')}.
\end{aligned}
\end{equation}
Newton's rule was originally given in terms of densities,  as above. The formulation (\ref{eq:newton}) is, however, more convenient for our purposes.

\cite{newtonQuintanaZhang1998} first propose the recursive rule in the context of interval censored data and mixtures of Markov chains, further developed by \cite{newtonZhang1999}. Theoretical properties have been studied from a frequentist viewpoint,  that is, regarding $G_n$ as an estimator of the mixing distribution under the assumption that the data are i.i.d. according to a true (identifiable) mixture model.  \cite{smithMakov1978} prove frequentist consistency of their recursive estimator for finite mixtures, using stochastic approximation techniques. 
\cite{martinGhosh2008} shed light on the connection with stochastic approximation, thus relating frequentist consistency of the algorithm to the convergence properties of stochastic approximation sequences. 
%They prove consistency of $G_n$ for a discrete $G_{true}$ with known atoms, and extend to the case of mixture kernels with an unknown common parameter. 
\cite{ghoshTokdar2006} and \cite{tokdarMartinGhosh2009} prove frequentist weak consistency of the estimator (\ref{eq:newton}) under conditions on the mixture kernels.
% and give results on convergence in probability for a permutation-invariant version of it. 
These results regard Newton's algorithm (\ref{eq:newton}) as a frequentist estimator. Its properties as an approximation of a  computationally expensive Bayesian solution remain unknown.  
One could argue that, when consistent, Newtons' estimator will asymptotically agree, almost surely with respect to the law of i.i.d. observations from a true $F_{G_{true}}$, with a consistent Bayesian estimator for $\tG$. But this is also true for any other consistent estimator of $\tG$, making them indistinguishable under this criterion. Newton's recursive estimator has the advantage of being  computationally faster than other consistent estimators for the mixing distribution, but its Bayesian motivation is lost. 
 
We take a different approach, focusing on its properties as a quasi-Bayesian procedure. 
First, we notice that the mixture model (\ref{eq:mixture}) can be expressed in terms of a latent sequence of random variables $(\theta_i)$ such that, given $(\theta_i)$, the $X_i$ are independent, with $X_i \mid \theta_i \sim f(x \mid \theta_i)$ and the $\theta_i$ are a random sample from $\tilde G$, i.e., $\theta_i \mid \tilde G \iidsim \tilde G$. Then, the Bayesian estimate $E(\tG \mid x_1, \ldots, x_n)$ coincides with the predictive distribution of $\theta_{n+1}$ given $(x_1, \ldots, x_n)$. 
Our point is that, when using (\ref{eq:newton}), a researcher is changing the predictive  rule of $\theta_{n+1}$,  therefore implicitly  using a  probabilistic model that is different from the Bayesian exchangeable model (\ref{eq:mixture}). What is this model? Is it quasi-Bayes? A similar reasoning and subsequent questions may arise in relation to other approximation algorithms which, more or less implicitly, use a probabilistic model different from the stated Bayesian one. 

%To address these questions, we need to 
Let us first formalize a notion of  {\em quasi-Bayes} procedure.  
The term {\em quasi-Bayes} is given many meanings in the literature (see e.g. \cite{liGuedjLoustau2018}). We borrow this term from \cite{smithMakov1978}, and formalize its meaning as follows.
A predictive rule implicitly defines the probability law, say $P$, of the sequence $(X_n)$. 
%The basic setting for Bayesian inference is exchangeability of the infinite sequence $(X_n)$. A quasi-Bayes learning rule should preserve this invariance property, at least asymptotically. That is, it should be {\em asymptotically exchangeable}. Informally, the observations' labels do not matter if we look at sufficiently distant pieces of the sequence. 
%Yet, asymptotic exchangeability is only a minimal property. A refinement refers to a specific exchangeable law. If $\tP$ is an exchangeable probability law for $(X_n)$, 
We say that $P$ is a quasi-Bayes approximation of an exchangeable probability law $\tP$  if it is asymptotically exchangeable,  and the exchangeable limit sequence has probability law $\tP$.  On this basis,  we address the following questions: 

1. If  one uses (\ref{eq:newton}) as a probabilistic learning rule, that is, as the predictive distribution of $\theta_{n+1}$ given $(X_1, \ldots, X_n)$, what statistical model is she implicitly assuming for the observable $(X_n)$?  Is it an approximation, at least asymptotically, of a Bayesian, exchangeable, mixture model? 

2. As an algorithm, Newton's recursive rule (\ref{eq:newton}) provides only a point estimate of the mixing distribution %$G$. 
Instead, a quasi-Bayes method should fully describe the uncertainty through the posterior distribution. Can  a posterior distribution %on $G$ 
be provided?

It is well known that the recursive estimate $G_n$ is not invariant to permutations of the observations. This means that the underlying probability law is not exchangeable. In fact, we  show that it implies a weaker form of dependence; namely, 
%We prove that the probabilistic model underlying the recursive rule (\ref{eq:newton}) is asymptotically exchangeable. More precisely, it implies that 
the sequence $(X_n)$ is {\em conditionally identically distributed} (c.i.d.; \cite{Kallenberg88}, \cite{BertiPratelliRigo2004LimThForIID}). Roughly speaking, for any $n$, future observations $X_{n+k}, k \geq 1$, are identically distributed, given $(X_1, \ldots, X_n)$ (see Section \ref{sec:pred-cid}). 
 For stationary sequences, the c.i.d. property is equivalent to exchangeability, while in general, a c.i.d. sequence is only asymptotically exchangeable. Therefore, a researcher using (\ref{eq:newton}) as the predictive rule is implicitly assuming some form of non-stationarity in the data, which tends to vanish  in the long run. A c.i.d. model could be the appropriate model in situations where exchangeability is broken by competition, selection (or other forms of non stationarity), but the system converges to a stationary, exchangeable, steady state. We develop this notion in a time-dependent mixture model in Section \ref{sec:estensioni}.  
 If, instead, the c.i.d. model is used as an approximation of an exchangeable model, it guarantees the minimal property of being asymptotically exchangeable: informally, 
  for $n$ large, the law of $ (X_{n+1}, X_{n+2}, \ldots)$ is approximately invariant under permutations, and there exists a random distribution, say $\tF$, such that $X_i \mid \tF \approxiid \tF$ for $i>n$, where $\approxiid$ denotes approximately i.i.d. Such $\tF$  plays the role of the asymptotic statistical model. % for the sequence.  

We refine this result by finding the explicit form of the asymptotic statistical model $\tF$ implied by the recursive predictive rule. The lack of exchangeability of the sequence $(\theta_n)$ implies that there is no random distribution $G$ such that $\theta_i \mid G \iidsim G$; however, we prove that such $G$ exists asymptotically, as the almost sure weak limit of the sequence $G_n$. Then the asymptotic statistical model 
%We prove that, almost surely, there exists a random distribution $G$, which is 
%the almost sure weak limit of the sequence $G_n$, such that 
$\tF$ is a mixture of the form $F_G=\int f(x \mid \theta) dG(\theta)$; roughly speaking, $X_i \mid G \approxiid F_G$ for $i > n$ with $n$ large.
In  this sense, Newton's recursive learning rule arises from a {\em quasi-Bayes} mixture model. 

At first, this result may appear surprising, as we are reverting the usual approach that goes from a statistical model and a prior distribution to the consequent predictive distribution. Here, we start from the predictive distribution (the algorithm) and find the implied statistical model. The key is that c.i.d. sequences preserve the main asymptotic properties of exchangeable sequences; in particular, the asymptotic statistical model arises as the almost sure weak limit of the predictive distribution. % of $X_{n+1}$ given $X_{1:n}$. }
%The key is that, for exchangeable sequences $(X_n)$, the predictive distribution almost surely converges weakly to a random distribution, and by de Finetti's representation theorem, the $X_n$ are conditionally i.i.d. according to such limiting distribution. Thus, one can find the statistical model as the random limit of  the predictive distributions, and the prior as the probability law of such random limit. A similar property holds for c.i.d. sequences: the predictive distribution of $X_{n+1}$ given $X_{1:n}$ converges to a random distribution, which acts as the {\em asymptotic} statistical model. 
The proofs of our results (collected in the Appendix) %Section \ref{sec:proofs})
come  from this property. In fact, we show a stronger result: the predictive distribution $F_{G_n}$ converges almost surely in total variation to the mixture $F_G$ (See Section \ref{sec:quasi-Bayes}). 

These results shed light on an open question posed by \cite{martinGhosh2008}. Although Newton's algorithm is popularly used for approximate computations in  DP mixture models, the authors show two examples where the Bayesian estimate of $G$ with a DP prior and the recursive estimate $G_n$ have different performance. Thus, they pose question:  {\em If Newton's recursive algorithm is not an approximation of the DP prior Bayes estimate, for what prior does the recursive estimate approximate the corresponding Bayes estimate?} 
We have shown that the latent distribution $G$ exists only asymptotically, so that the prior has to be interpreted as the probability law of such {\em asymptotic} $G$; and will generally differ from a DP. 
%We explain the reason why the prior law may differ from a DP.  The recursive rule (\ref{eq:newton}) underlies a sequence $(\theta_n)$ that is not exchangeable; thus, there is no random mixing distribution $G$ such that $\theta_i \mid G \iidsim G$. Such $G$ only exists asymptotically, as the weak limit of the predictive distributions $G_n$. 
%and its probability law will generally differ from a DP.
While the DP is almost surely discrete, we prove that,  under fairly mild conditions, the random distribution $G$ is absolutely continuous, almost surely. Thus, the c.i.d. model implies a novel prior on absolutely continuous random distributions.
%: informally, a smoothing of the Dirichlet process. 

This clean statistical formulation of Newton's algorithm allows to provide a probabilistic description of uncertainty, through the posterior distribution.  A different proposal, %way to express uncertainty for Newton's algorithm, 
based on the variability of the mixing density estimates obtained over random permutations of the original data, has been recently suggested by \cite{dixitMartin2019}. 
We aim at a proper posterior distribution. 
Indeed,  although the prior law of $G$ is only implicitly defined, we can approximate the corresponding posterior distribution of $G$, leveraging on  properties of c.i.d. processes. 
 More precisely, we obtain an asymptotic Gaussian approximation of the posterior distribution $P([G(A_1), \ldots, G(A_k)] \in \cdot \mid x_1, \ldots, x_n)$ for any measurable sets $A_1, \ldots, A_k$, $k \geq 1$. 
Thus, additionally to a quasi-Bayes point estimate, one may provide asymptotic credible regions. 

These results develop Newton's algorithm into a quasi-Bayes statistical method for sequential learning in mixture models, and shed light on 
the role of the different components of the model. In particular, we show that the weights $\alpha_n$ in the recursions have a dual role as learning parameters in the predictive rule and  as parameters that regulate the speed of convergence to approximate exchangeability. Consequently, they not only affect the sensitivity of the estimates to permutations of the data, but also the width of the  asymptotic credible intervals.  Through synthetic examples, we discuss this trade off in the choice of the weights and suggest  practical hints for tuning them. 
Interestingly, as we discuss in Section \ref{sec:inferenceG}, considering the balance between the learning rate and the predictive convergence rate provides novel insights on frequentist coverage, that apply, more generally, to the exchangeable setting. 

Another relevant implication of our approach %statistical approach for interpreting the algorithm 
is that we can naturally envisage extensions in several directions. The original version of Newton's algorithm does not cover  the case of common unknown parameters in the mixture components. Extensions are proposed  by \cite{martinGhosh2008} for  finite mixtures with a known support of the mixing distribution. The procedure suggested by \cite{martinTokdar2011} is  more general, but remains somehow heuristic, due to the lack of a genuine likelihood.
 Our results provide such likelihood and allow to properly envisage empirical Bayes or Bayesian inference of common unknown parameters. Moreover, we can pursue quasi-Bayes inference on the individual parameters $\theta_i$. Our results show that, in Newton's model, the latent distribution may be $P$-a.s. absolutely continuous, with density $g$. In this case, one obtains
multiple-shrinkage effects in the estimation of the $\theta_i$, guided by the modes of the latent density $g$.
A known limitation of Newton recursions is that they require to evaluate the normalizing constant at each step. While numerical integration is effective for a low-dimensional parameter $\theta$, it become cumbersome in the multivariate case.  We suggest a  recursive Monte Carlo sampling scheme to overcome this difficulty. Preliminary results we obtain in  simulation studies are encouraging. 

\medskip

We set the notation and remind preliminary notions in Section \ref{sec:DPNewton}. 
%Predictive constructions and c.i.d. sequences are recalled in section ...... 
%on Dirichlet process mixture models and predictive constructions in Section \ref{sec:DPNewton}. 
The proposed predictive methodology is developed in  Section \ref{sec:main}, where we provide a statistical interpretation of Newton's algorithm and find the implied modeling assumptions.  
These results are used in Section 4 to obtain an asymptotic approximation of the posterior distribution of the mixing distribution and the corresponding credible intervals. 
In Section \ref{sec:estensioni} we define a time-dependent mixture model consistent with the recursive predictive rule and discuss the role of the model parameters via simulation studies. We provide further statistical applications in  Section \ref{sec:moreExamples}. In Section \ref{sec:discussion} we briefly discuss future lines of research. All the proofs are collected in the Appendix. 

\section{Preliminaries: Dirichlet process mixtures and predictive characterizations}
    \label{sec:DPNewton}
    
Let us first set some notation. 
We consider random variables $X_i \in \mathbb X \subseteq \mathbb{R}^d$ and 
$\theta_i \in \Theta \subseteq \mathbb{R}^p$ (but our results hold for general  Polish spaces), where $\mathbb X$ and $\Theta$ are 
 equipped with the Borel sigma-fields $\mathcal B(\mathbb  X)$ and $\mathcal B(\Theta)$.
Throughout the paper, we refer to conditional distributions as regular versions. We use the short notation $X_{1:n}=(X_1, \ldots, X_n)$, and $P(A \mid x_{1:n})$ for $P(A \mid X_1=x_1, \ldots, X_n=x_n)$. A sequence $(Z_i)_{i=1}^n$ will be briefly written as $(Z_n)$.
We use the same symbol to denote a probability measure and the corresponding distribution function. Unless explicitly stated, weak convergence of distributions is considered and denoted by $F_n\Rightarrow F$.
    
%Let us first recall 
We now briefly remind the basic structure of Bayesian inference for DP mixture models, 
in order to motivate in more detail the recursive rule (\ref{eq:newton}).
%and to introduce some further notation. 
Again, a DP mixture model has a hierarchical formulation in terms of a latent exchangeable sequence $(\theta_i)$
\begin{eqnarray} 
&& X_i \mid \theta_i \indsim f(x \mid \theta_i)  \label{eq:short}\\
&& \nonumber \\ 
&& \theta_i \mid \tG \iidsim \tG, \nonumber 
\end{eqnarray}   
with $\tG \sim \DP(\alpha G_0)$,  
where (\ref{eq:short}) is abbreviated notation for 
$X_n \mid X_{1:n-1}, (\theta_n) \sim f(x\mid\theta_n)$, for every $n \geq 1$, and $f(\cdot \mid \theta)$ is a density with respect to a sigma-finite measure $\mu$ on the sample space $\mathbb X$. Integrating the $\theta_i$ out, one has the mixture model (\ref{eq:mixture}), with a  \DP prior on $\tG$.  We denote by $\tP$ the probability law on the process $( (X_n, \theta_n))_{n \geq 1}$ so defined. 
%%%%
%We assume that the $X_i$ and the $\theta_i$ are random variables with values in $\mathbb X \subseteq \mathbb{R}^d$ and $\Theta \subseteq \mathbb{R}^p$, respectively, equipped with the Borel sigma-fields $\mathcal B(\mathbb  X)$ and $\mathcal B(\Theta)$;  but the results hold for general  Polish spaces. We use the same symbol to denote a probability measure and the corresponding distribution function. Throughout the paper, we refer to conditional distributions as regular versions. We use the short notation $X_{1:n}=(X_1, \ldots, X_n)$, and $P(A \mid x_{1:n})$ for $P(A \mid X_1=x_1, \ldots, X_n=x_n)$. 
%Unless explicitly stated, weak convergence of distributions is considered and denoted by $F_n\Rightarrow F$.
%%%%
Inference on the latent distribution $\tG$ in a \DP mixture model is solved through the posterior distribution, which is a mixture of DPs \citep{antoniak1974}
\begin{equation} \label{eq:posteriorDPmixture}
\tG \mid x_{1:n} \sim \int \DP(\alpha G_0 + \sum_{i=1}^n \delta_{\theta_i}) \,  d \tP(\theta_{1:n} \mid x_{1:n}). 
\end{equation}
The Bayesian point estimate $\tG_n^{(Bayes)}$ of $\tG$, with respect to quadratic loss,  is the conditional expectation of $\tG$, and coincides with the predictive distribution of $\theta_{n+1}$, given $X_{1:n}$. 
By the \Polya urn structure characterizing the Dirichlet process
\begin{equation} \label{eq:PolyaDP}     
\tP(\theta_{n+1} \in \cdot \mid \theta_{1:n}, x_{1:n}) = \frac{\alpha G_0(\cdot) + \sum_{i=1}^n \delta_{\theta_i}(\cdot)}{\alpha+n}, 
\end{equation}
therefore
\begin{eqnarray} 
\tG_n^{(Bayes)}(\cdot)&=&E\tG(\cdot) \mid x_{1:n})= \tP(\theta_{n+1} \in \cdot \mid x_{1:n}) = \frac{\alpha G_0(\cdot) + \sum_{i=1}^n \tP(\theta_i \in \cdot \mid x_{1:n})}{\alpha+n}  \label{eq:PolyaMixDP} \\
&=& \frac{\alpha+n-1}{\alpha+n} \, 
\frac{\alpha G_0(\cdot) + \sum_{i=1}^{n-1} \tP(\theta_i \in \cdot \mid x_{1:n})}{\alpha+n-1} + \frac{1}{\alpha+n}   \tP_{\tG_{n-1}^{Bayes}}(\theta_n \in \cdot \mid x_n), \nonumber    
\end{eqnarray}
where  we use the notation
\begin{equation} \label{eq:posteriorDP}
\tP_{G}(\theta_n \in A \mid x_n) = 
\frac{\int_A f(x_n \mid \theta) d G(\theta)}{ \int_\Theta f(x_n \mid \theta) d G(\theta)}. 
\end{equation}
In the Bayesian estimate, as a new observation $x_n$ becomes available, the information on all the past $\theta_i, i=1,\ldots, n-1$, is updated. This efficiently exploits the sample information, but is computationally expensive. Instead, Newton's algorithm (\ref{eq:newton}) does not update the estimate $G_{n-1}$, and $x_n$ only enters the inference on $\theta_n$, with an empirical Bayes flavor. The two estimates coincide only for $n=1$ and, even in this case,
Newton's rule makes a simplification of the posterior distribution of $G$, replacing the mixture of Dirichlet processes  $\int \DP(\alpha G_0+\delta_{\theta_1}) d\tP(\theta_1 \mid x_1)$, as from  (\ref{eq:posteriorDPmixture}),  with a \DP($\alpha G_0 + \tP(\theta_1 \mid x_1))$. For $n \geq 1$, Newton's estimate loses  efficiency, not fully exploiting the sample information. On the other hand, it is very fast. If one evaluates (\ref{eq:newton}) on a grid of $m$ points and calculates the integral in the denominator using, say, a trapezoid rule, then the computational complexity is $m n$. 
%\red{(da Martin e Ghosh, pag 373)}

\subsection{Predictive constructions and conditionally identically distributed sequences}
\label{sec:pred-cid}

As anticipated, the key of our developments is to regard the recursive rule (\ref{eq:newton}) as a probabilistic {\em predictive rule}.  Let us briefly remind the essentials of the predictive approach to inference highlighting an interesting form of stochastic dependence that emerges from it, namely the notion of {\em conditionally identically distributed} sequences. 

In Bayesian inference, predictive characterizations are a natural and powerful 	tool to define prior distributions. We mention the well known predictive characterization of the DP through \Polya sequences \citep{blackwellMacqueen1973}, or Chinese Restaurant Process, and we refer to \cite{fortiniP2012b} for a review. Let $(Z_n)$ be a sequence of random variables, and for any $n \geq 1$ let $P_n(\cdot)=P(Z_{n+1} \in \cdot \mid Z_{1:n})$. By the Ionescu-Tulcea theorem, the sequence $(P_n)$ characterizes the probability law, $P$, of $(Z_n)$. If $P$ is exchangeable, then, by de Finetti representation theorem, it characterizes the implied prior law. An interesting result by \cite{Kallenberg88} (Proposition 2.1) proves that a stationary sequence satisfying
\begin{equation} \label{eq:cidKallenberg}
 (X_1, \ldots, X_n, X_{n+2}) \stackrel{d}{=}(X_1,\dots,X_n,X_{n+1}), \quad n \geq 1,
\end{equation}
where $\eqdist$ means equal in distribution, is exchangeable. Clearly, the converse is true, thus condition (\ref{eq:cidKallenberg}) is equivalent to exchangeability for stationary sequences. Therefore, a predictive rule characterizes an exchangeable probability law $P$ if and only if $P$ is stationary and satisfies (\ref{eq:cidKallenberg}). 

Notice that (\ref{eq:cidKallenberg}) implies that $(X_1, \ldots, X_n, X_{n+k}) \stackrel{d}{=}(X_1,\dots,X_n,X_{n+1})$, for any $n \geq 1$ and $k\geq 1$. Informally, for any $n \geq 1$, 
$$X_{n+k}\mid X_{1:n} \eqdist X_{n+1}\mid X_{1:n}, \quad \mbox{for any $k \geq 1$}. 
$$ 
  \cite{BertiPratelliRigo2004LimThForIID} extend this notion, introducing the notion of {\em conditionally identically distributed sequences with respect to a filtration} and provide fundamental limit theorems. 
Interestingly, c.i.d. sequences preserve the main asymptotic properties of exchangeable sequences. In particular,  the sequence of the empirical distributions and the sequence of the predictive distributions  converge $P$-a.s., to the same random distribution, i.e., if $(X_n)$ is c.i.d. with probability law $P$, then 
\begin{equation} \label{eq:convExch}
\hat{F}_n \equiv 	\frac{\sum_{i=1}^n \delta_{X_i}}{n} \Rightarrow F \; \;  \mbox{and} \; 
P_n \equiv P(X_{n+1} \in \cdot \mid X_{1:n}) \Rightarrow F , \; \; \mbox{$P$-a.s.}.
\end{equation} 
For exchangeable sequences, the limit $F$ is called the directing random measure  (the statistical model, in Bayesian inference) and the probability law of $F$ is the de Finetti measure (the prior distribution). The term
directing random measure is used analogously for c.i.d. sequences.

An exchangeable sequence is clearly c.i.d., but the converse is not generally true. However, c.i.d. sequences are asymptotically exchangeable. 

\begin{defn}  
A sequence of random variables $(X_n)$ is asymptotically exchangeable, with directing random measure $F$, if  
$$ (X_{n+1}, X_{n+2}, \ldots) \overset{d}{\rightarrow} (Z_1, Z_2, \ldots) $$  
for an exchangeable sequence $(Z_n)$, with directing random  measure $F$.  
\end{defn}

For a sequence $(X_n)$, convergence of the predictive distributions to a random probability measure, $\mu$,  implies that the sequence is asymptotically exchangeable with directing random measure $\mu$
(\cite{aldous1985} Lemma 8.2).  
Thus, by (\ref{eq:convExch}), a c.i.d.  sequence $(X_n)$ is asymptotically exchangeable, with directing random measure $F$. Informally, $X_n \mid F {\overset{ind}\approx} F$, for large $n$. 

%Fundamental results and limit theorems for $\cF$-c.i.d. sequences are given by \cite{BertiPratelliRigo2004LimThForIID}. 
Applications of c.i.d. processes in Bayesian nonparametric inference  include \cite{Bassetti2010CidSSS} and the c.i.d. hierarchical model proposed by \cite{airoldiCosta2014}. 
 
\section{A statistical interpretation of Newton's algorithm} \label{sec:main}  

For the mixture model (\ref{eq:short}), the Bayesian point estimate $E(\tG \mid x_{1:n})$ corresponds to the predictive distribution of $\theta_{n+1}$ given $x_{1:n}$. Our point is that, similarly, Newton's recursive rule (\ref{eq:newton}) should be regarded as a different {\em probabilistic predictive distribution} for $\theta_{n+1}$ in the latent variable model (\ref{eq:short}), assuming  
\begin{equation} \label{eq:newtonAsPred}
\begin{array}{cc}
X_n \mid \theta_n \indsim f(x \mid \theta_n) \\
\theta_{n+1} \mid x_{1:n} \sim G_n(\cdot),  \; n \geq 1,
\end{array}
 \end{equation}   
with $\theta_1 \sim G_0$ and $G_n$ given by (\ref{eq:newton}). According to the predictive approach (Section \ref{sec:pred-cid}), this means that a researcher using the recursive rule (\ref{eq:newton}) is implicitly assuming  a different statistical model for the sequence $((X_n, \theta_n))_{n \geq 1}$, in place of the exchangeable mixture model (\ref{eq:mixture}) and it is important to make such a  model explicit. This model may be of autonomous interest in some experimental circumstances. 
%When regarded as an approximation of a Bayesian procedure, it provides the formal framework for understanding the quasi-Bayes properties of the recursive rule.

Let us denote by $P$ a probability law on the joint process $((X_n, \theta_n))_{n \geq 1}$  that is consistent with the assumptions (\ref{eq:newtonAsPred}). 
A first implication of our approach is that the recursive formulae can now be given a probabilistic interpretation. The estimate $G_n$ can be written in a prediction-error correction form,
$$
G_{n}(\cdot)=G_{n-1}(\cdot)+\alpha_n [P(\theta_{n}\in \cdot \mid x_{1:n})-P(\theta_{n}\in\cdot\mid x_{1:n-1})], 
$$
where the correction term is now properly interpreted as a difference between predictive distributions computed according to $P$. Moreover, we can appreciate the different information conveyed by the recursive predictive rule with respect to DP mixtures. Simple computations show that one can write $G_n$ as  
\begin{equation} \label{eq:defGn}
G_n(\cdot)=\frac{\alpha G_0(\cdot)+\sum_{k=1}^n\gamma_k P(\theta_k\in\cdot\mid x_{1:k})}{\alpha+\sum_{k=1}^n\gamma_k}, \quad \quad n\geq 1, 
\end{equation}
where $P(\theta_k\in\cdot\mid x_{1:k})= P_{G_{k-1}}(\cdot \mid x_k)$;  $\alpha>0$, $\gamma_1=\alpha_1\alpha/(1-\alpha_1)$ and $\gamma_n=\alpha_n(\alpha+\sum_{k=1}^{n-1}\gamma_k)/(1-\alpha_n)$ for $n\geq 2$. 
For $\alpha_n = 1/(\alpha+n)$, one has $\gamma_n=1$ for all $n \geq 1$ and a direct comparison with the corresponding formula (\ref{eq:PolyaMixDP}) for DP mixtures.  

The rule (\ref{eq:PolyaMixDP}) originates from the \Polya urn scheme characterizing the Dirichlet process. This suggests that Newton's recursions are based on a different urn scheme, possibly an urn of distributions (see \cite{quintanaNewton2000}). Having framed the recursive rule in a probabilistic setting, we can make such intuition rigorous, by providing the explicit form of the predictive rule for the observable sequence $(X_n)$. Interestingly, it proves be a novel {\em measure-valued \Polya urn scheme} (\cite{bandyopadhyaThacker2017}, \cite{maillerMarkert2017}, \cite{janson2019}).
From (\ref{eq:newtonAsPred}) it follows that $X_1 \sim F_{G_0}(\cdot) \equiv  \int F(\cdot \mid \theta) dG_0(\theta)$ and for any $n \geq 1$
$$ X_{n+1} \mid x_{1:n} \sim P_n (\cdot) = \int F(\cdot \mid \theta) dG_n(\theta) = (1-\alpha_n) P_{n-1}(\cdot) + \alpha_n F_{G_{n-1}}(\cdot \mid x_n) , $$
where $F(\cdot\mid\theta)$ is the distribution function with density $f(x\mid\theta)$ with respect to $\mu$ and $F_{G_{n-1}}(\cdot \mid x_n)= \int F(\cdot \mid \theta) d P_{G_{n-1}}(\theta \mid x_n)$. 
This novel \Polya urn scheme provides a predictive characterization of the probability law of the process $(X_n)$.  When using Newton's rule, a researcher should be aware of the  assumptions made on the observable $(X_n)$ through such a probability law, and making such assumptions explicit is the aim of the next section. 

\subsection{Quasi-Bayes properties} \label{sec:quasi-Bayes}
 
Newton's model (\ref{eq:newtonAsPred}) does not fully specify the probability law of the process $( (X_n, \theta_n) )$,  because it only assigns the probability law of $\theta_{n+1} \mid X_{1:n}$ and not enough restrictions are made on the conditional distributions of $\theta_{n+1} \mid X_{1:n}, \theta_{1:n}$. Nevertheless, it has relevant implications, which we study in this section. 
Clearly,  a trivial way to obtain a full specification is to assume that $\theta_{n+1}$ is conditionally independent on  $\theta_{1:n}$, given $X_{1:n}$. This might be motivated by the non-stationary nature of the sequence $(\theta_n)$ in (\ref{eq:newtonAsPred}), and would simplify the analysis but is an unnecessary additional assumption, as our results show.

As noticed in the Introduction, the probability law implied by the assumptions (\ref{eq:newtonAsPred}) is not exchangeable;  we  show that the process $(X_n)$ is in fact c.i.d. 
Yet, a mixture model of the form $X_n \mid G \iidsim f_G$, at least  asymptotically, is desirable. 
%It is well known that Newton's recursive estimate is not invariant under permutations of the observations. In the probabilistic setting, this means that the probability law implied by the assumptions (\ref{eq:newtonAsPred}) is not exchangeable, we  show that the process $(X_n)$ is in fact c.i.d. In some applications, exchangeability may be broken by forms of non-stationarity, and a c.i.d. model may offer a sensible description of the phenomena. We  return on this case in Section \ref{sec:timeMixture}.  In a standard setting, the researcher would regard the $X_i$ as exchangeable, judging that the observations' labels do not matter. The lack of exchangeability implied by assumptions (\ref{eq:newtonAsPred}) is thus a misspecification, motivated only by the need of fast computations. Then  $(X_n)$ should be at least asymptotically exchangeable: informally, the labels do not matter if we look at $(X_{n}, X_{n+1}, \ldots)$, for  large enough $n$. Moreover, a mixture model of the form $X_n \mid G \iidsim f_G$, at least  asymptotically, is desirable. To prove that the latter holds, 
We start by showing that an asymptotic mixing distribution $G$ exists, and is the $P$-a.s. limit of the sequence of the predictive distributions $G_n$. Furthermore, 
% as for exchangeable sequences, 
$G_n$ is the conditional expectation of $G$, given $X_{1:n}$. The proofs of the theorem below and  of all subsequent results are collected in the Appendix.
%Section \ref{sec:proofs}.  

\begin{thm} \label{th:urnmart}
Let the process $((X_n, \theta_n))$ have a probability law $P$ that satisfies assumptions (\ref{eq:newtonAsPred}). Then, $P$-a.s.,  
\begin{itemize}
\item[(i)]	the sequence $(G_n)$ converges to a random probability measure $G$,
\item[(ii)] for every $n \geq 1$ and measurable set $A$, 
$P(\theta_{n+k} \in A \mid X_{1:n})= E( G(A) \mid X_{1:n})$, for all $k \geq 1$. 
\end{itemize}
\end{thm}
 
An immediate consequence of the weak convergence of $G_n$ to $G$ is that $ \int h(\theta)dG_n(\theta) \rightarrow \int h(\theta) dG(\theta)$ $P$-a.s. for any continuous and bounded function $h$ on $\Theta$. We prove  that the convergence can be extended to functions $h$ that are  integrable with respect to $G$. 

\begin{prop}	\label{prop:hconv}
Let $((X_n, \theta_n)) \sim P$ satisfy the  assumptions (\ref{eq:newtonAsPred}), and let $h(\cdot)$ be a measurable function on $\Theta$, such that,  $P$-a.s., 
$\int |h(\theta)| dG(\theta) < \infty$. 
Then, for $n \rightarrow \infty$, 
$$ 
\int h(\theta) dG_n(\theta) 
\rightarrow \int h(\theta) dG(\theta), \quad \mbox{$P$-a.s.}
$$
The condition $\int |h(\theta)| dG(\theta) < \infty$ $P$-a.s. holds, in particular, if $h$ is measurable and $\int |h(\theta)| dG_0(\theta) < \infty$.
\end{prop}
 
The following theorem proves that Newton's learning rule (\ref{eq:newtonAsPred}) implies that the sequence $(X_n)$ is c.i.d., thus asymptotically exchangeable, and that its directing random measure has a mixture density of the form $f_G$. 
In this sense, Newton's model is a {\em quasi-Bayes} mixture model. 

\begin{thm}  	\label{th:cidconv}
Let $((X_n, \theta_n)) \sim P$ satisfy  the assumptions (\ref{eq:newtonAsPred}). 
	Then 
\begin{itemize}
\item[(i)] The sequence ($X_n$) is c.i.d.; 
\item[(ii)]  The sequence of predictive densities $f_{G_n}$ converges in $L_1$ to 
$f_G \equiv \int f(x \mid \theta) dG(\theta)$, $P$-a.s.,  where $G$ is the $P$-a.s. weak limit of $(G_n)$; 
 \item[(iii)] 
$(X_n)$ is asymptotically exchangeable, and its directing random measure has density $f_G$ with respect to $\mu$.
\end{itemize}
\end{thm}
 
Informally, the above results say that $X_n \mid G \approxiid f_G$, for $n$ large. Notice that $G$ plays the role of the (infinite-dimensional) parameter of the asymptotic statistical model $F_G$ of $(X_n)$, and,  as such, it is a function of $(X_n)$. If the mixture is identifiable, then $F_G$ uniquely determines $G$. 
Moreover, by properties of c.i.d. sequences, $F_G$ is also the $P$-a.s. weak limit of the sequence of empirical distributions $\sum_{i=1}^n \delta_{X_i}/n$. 
%Therefore, it is measurable with respect to the tail sigma-field of $(X_n)$. 

Intuitively, asymptotic exchangeability of the sequence $(X_n)$ implies that  also $(\theta_n)$ is asymptotically exchangeable. In fact, if we assume the additional condition that $\theta_{n+1}$ is independent on $\theta_{1:n}$ given  $X_{1:n}$, then  it is easy to prove that $(\theta_n)$ is c.i.d., thus asymptotically exchangeable, but such assumption is not necessary, as shown by the following theorem. 
 
\begin{thm}
	\label{th:asympttheta}
If the mixture $F_G=\int f(y \mid \theta) dG(\theta)$ is identifiable, then Newton's learning scheme (\ref{eq:newtonAsPred}) implies that the sequence $(\theta_n)$ is asymptotically exchangeable, with directing random measure $G$, corresponding to the $P$-a.s. limit of the sequence $G_n$.
\end{thm}
 
\subsection{On the prior distribution of $G$}  \label{sec:prior}

%While arising from the first step of the updating rule in DP mixtures,  which assume an a.s. discrete mixing distribution, Newton's algorithm was proposed as a recursive estimate of the mixing {\em density}, implicitly assuming that the unknown  mixing distribution $G$ is absolutely continuous.  
We  have shown that, when interpreted as predictive rules, Newton's recursions imply a quasi-Bayes mixture model. Yet, in contrast with DP mixtures, the prior on $G$ is no longer, in general, a DP. 
Explicit results on the probability law of the random limit distribution $G$ are challenging. In the probabilistic literature on c.i.d. processes, very few results of this nature are available, for very simple cases. Even in the exchangeable case, finding the prior implicitly defined through a predictive construction is often an open problem. 
Although we cannot provide the explicit form of  the prior distribution on $G$ implied by Newton's model (\ref{eq:newtonAsPred}), we can prove that, under fairly mild sufficient conditions, $G$ is $P$-a.s. absolutely continuous. 
Moreover, in the next section we provide an asymptotic Gaussian approximation of the posterior distribution that results from such implicit prior.   

If $G_0$ is absolutely continuous with respect to a sigma-finite measure $\lambda$ on $\Theta$,  denoted $G_0 \ll \lambda$, then $G_n \ll \lambda$, and the corresponding density $g_n$ satisfies Newton's recursive rule (\ref{eq:recursiveDensities}). 
It is easy to verify that, for any fixed $\theta$,  the sequence $(g_n(\theta))$ is a martingale under the c.i.d. law $P$.  Since $g_n(\theta)$ is non-negative, there exists a function $g^*(\theta)$ such that, for every $\theta$,  $g_n(\theta)$ converges to $g^*(\theta)$,  $P$-a.s. However, this fact is not sufficient to conclude  that $G \ll \lambda$. 
%Extending a remarkable result by \cite{BertiPratelliRigo2013}, we show that $G \ll \lambda$ requires that $g_n$ converges in $L^1$ or, equivalently, that $G_n$ converges to $G$ in total variation.  
%This is guaranteed in Lemma \ref{th:appendix1}. 
%Then, Lemma \ref{th:appendix2} gives sufficient conditions for $G_n \rightarrow G$ in total variation. % and $G\ll\lambda$, a.s.
%The two lemmas are extensions of Theorems 1 and 4 in  \cite{BertiPratelliRigo2013}, the main difference being that we do not assume that the sequence of random variables involved is c.i.d.  
%\red{(a property that does not hold, in general, for the sequence $(\theta_n)$ in Newton's model).} 
%%\red{CUT? To be c.i.d. with respect to a filtration  $\cF$, a sequence has to be adapted to $\cF$, which, in general, is not true for $(\theta_n)$ in Newton's model and the filtration of interest.}
%However, the proofs of those theorems can be fairly easily adapted  by directly requiring the martingale property of the sequence of random measures $Q_n$, which is otherwise implied by the $\cF$-c.i.d. property. 
%%However, the proofs of Theorems 1 and 4 in \cite{BertiPratelliRigo2013} can be adapted,
%% by directly requiring the martingale property of the sequence of the random measures $Q_n$, which is otherwise implied by the $\cF$-c.i.d. property.
%Thus, the following two lemmas are provided without additional proof.
The following theorem gives sufficient conditions for $G$ to be absolutely continuous with respect to $\lambda$. In particular, when $\lambda$ is the Lebesgue measure, it gives conditions for the existence of a density $g$, which turns out to be the limit of $(g_n)$, in $L^1$. 
	The proof is based on a remarkable result in \cite{BertiPratelliRigo2013} (Theorems 1 and 4), which  shows that, for c.i.d. processes, the directing measure is absolutely continuous with respect to $\lambda$ if and only if the predictive distribution is absolutely continuous and converges in total variation. 
However, this result does not  apply directly to our setting, since $(\theta_n)$ is not, generally, c.i.d. %such (it is not necessarily adapted to the filtration generated by $(X_n)$). 
In the Appendix 
%we extend Theorems 1 and 4 in \cite{BertiPratelliRigo2013} (Lemmas \ref{th:appendix1} and \ref{th:appendix2}), }
we provide two Lemmas (Lemmas \ref{th:appendix1} and \ref{th:appendix2}), which are slight extensions of Theorems 1 and 4 in \cite{BertiPratelliRigo2013}, 
the main difference being that we substitute the c.i.d. assumption with a martingale property that holds in our setting. 
%However, the proof of those theorems can be easily adapted, thus the Lemmas are provided without additional proof. 
These results lead to the following theorem.
% 
%We can now provide sufficient conditions for the a.s. absolutely continuity of the limit mixing distribution $G$ in Newton's model.
\begin{thm}
	\label{th:absolcont}
Let $G$ be the $P$-a.s. limit of the sequence of predictive rules $G_n$ defined by (\ref{eq:newton}), with $G_0 \ll \lambda$. If the following conditions hold
\begin{equation}\label{eq:suffcont1}    
\sum_n\alpha_n^2<\infty \; ,   \quad   \; 
\int_Kg_0(\theta)^2d\lambda(\theta)<\infty, \; \mbox{for every $K$  compact,}
\end{equation}
and 
\begin{equation} \label{eq:suffcont3}
\sup_{\theta_1,\theta_2\in K}\int \frac{f(x\mid \theta_1)^2}{f(x\mid\theta_2)}d\mu(x)<\infty ,   \; \; 
\mbox{ for every $K$ compact such that $\lambda(K)<\infty$}, 
\end{equation}
 then $G\ll\lambda$, $P$-a.s.
Moreover, $P$-a.s, $g_n$ converges in $L^1$ to $g \equiv dG/d\lambda$.
	\end{thm}
	
Assumptions (\ref{eq:suffcont1})  are quite natural. They hold, for example, if $\alpha_n=1/(\alpha+n)$ and $g_0$ is continuous or bounded. Assumption (\ref{eq:suffcont3}) is more delicate. It holds, for example, if $f$ is a Poisson density or a Gaussian density with fixed variance or a Gamma density with fixed shape parameter. %It does not hold, for example, if $f$ is a general Gaussian density.
 A similar assumption is considered in \cite{tokdarMartinGhosh2009}. 

\begin{figure}[h]
\centering 
%\begin{center}
\includegraphics[width=.7\textwidth, height=7cm]{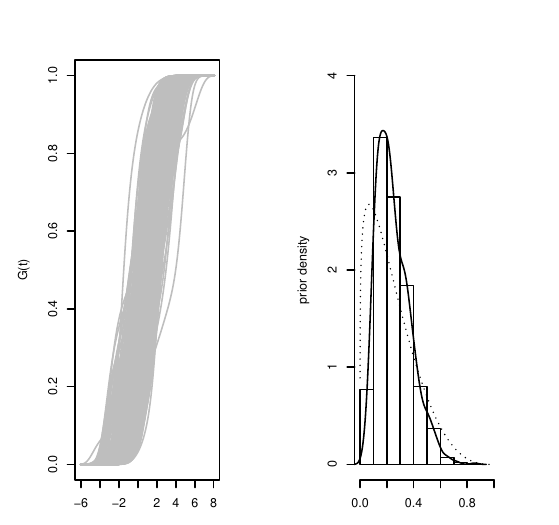}
%{Rplot-priorMC-new2}  % {plotPrior2-EX2.pdf}
% DATI SALVATI: 
% samplexS1 <- samplex
% sampleThetaS1 <- sampleTheta
% samplegnS1 <- samplegn
\caption{Monte Carlo approximation of the prior density of $G(0)$.
% versus the Beta prior density of $G(0)$ under the assumption that $G \sim \DP(\alpha G_0)$. 
First panel: Monte Carlo samples $G_N^{(m)}$, $m=1, \ldots, 1000$, $N=10,000$. 
%The Monte Carlo sample is obtained by sampling $M=1000$ replicates of data $x_{1:N}, N=10,000$, generated from the model (\ref{eq:newtonAsPred}) with a Gaussian kernel with known variance $\sigma^2=1$. The prior guess is $G_0= \Norm(1,2)$ and the weights are $\alpha_n=1/(\alpha+n)$ with $\alpha=5$. First panel: Monte Carlo samples $G^{(m)}_N$, $m=1, \ldots, 1000$. The vertical line shows the values $G^{(m)}(0)$. 
Second panel: Histogram of the sampled $G^{(m)}(0), m=1, \ldots, 1000$ and corresponding Monte Carlo estimate of the prior density of $G(0)$ (solid curve) versus the Beta$(\alpha G_0(0), \alpha (1-G_0(0))$ density (dotted).}
\label{fig:studioPrior}
\end{figure}

 \subsection{Empirical study} \label{sec:empirical-prior}
 
 Although we do not have the explicit expression of the prior law on $G$,  we can describe a procedure to simulate from it. Let us remind that the prior is the probability law of the random weak limit $G$ of the sequence $G_n$.   
For a continuity point $t$ of $G$,  one could, in principle, generate sequences  $\omega^{(m)}\equiv (x_1^{(m)}, x_2^{(m)}, \ldots)$ from the c.i.d. model, for $m=1, \ldots, M$, and for each $\omega^{(m)}$ compute $G^{(m)}(t) \equiv G(t)(\omega^{(m)})=\lim G_n(t)(\omega^{(m)})$. The resulting vector $[G^{(m)}(t), m=1, \ldots, M]$ would provide a Monte Carlo sample from the prior law of $G(t)$. Of course, one cannot generate an infinite sequence $(x_1^{(m)}, x_2^{(m)}, \ldots)$, but may use a truncated version $x_{1:N}^{(m)}$  with $N$ sufficiently large. 

In the following illustration, we consider Newton's model (\ref{eq:newtonAsPred}) with a Gaussian kernel $\Norm(\theta_i, \sigma^2)$ with known variance $\sigma^2=1$. 
The initial distribution $G_0$ is $\Norm(1,2)$ and the weights in the recursions are $\alpha_n=1/(\alpha+n)$ with $\alpha=5$. In this case, the assumptions of Theorem \ref{th:absolcont} are satisfied, and $G$ is $P$-a.s. absolutely continuous.  
%Figure \ref{fig:studioPrior}
We compare the Monte Carlo approximation of the prior density of $G(0)$ under Newton's model with the prior of $G(0)$ in a DP mixture model where $G \sim \DP(\alpha G_0)$, that is,
a Beta density with  parameters $(\alpha G_0(t), \alpha (1-G_0(t))$.
The Monte Carlo sample is obtained by generating $M=1000$ replicates $x_{1:N}^{(m)}$, $N=10,000$, from the c.i.d. model (\ref{eq:newtonAsPred}), assuming, for simplicity, that $\theta_{n+1}$ is independent on $\theta_{1:n}$ given $x_{1:n}$. For each sample $x_{1:N}^{(m)}$, we compute $G^{(m)}_N$, as a fairly reasonable proxy of the limit $G^{(m)}$. The first panel in Figure \ref{fig:studioPrior} shows the samples $G^{(m)}_N$ so obtained.  %emphasizing  their value at zero. 
The second panel shows the histogram of the sampled values $G^{(m)}_N(0)$ and a kernel density estimate, providing the Monte Carlo estimate of the prior density of $G(0)$. The dotted curve is the Beta$(\alpha G_0(0), \alpha (1-G_0(0))$ density that one would have under a DP prior on $G$. As expected, the two curves are quite different. 

%Indeed, by Theorem \ref{th:absolcont}, the prior under Newton's predictive rule selects a.s. absolutely continuous distributions under fairly mild conditions, while it is a.s. discrete under a DP prior. In statistical applications, DP mixture models are popular tools  for clustering, leveraging the implied random partition of the $\theta_i$. In Newton's model, the absolutely continuous $G$, instead determines multiple shrinkage effects in the estimation of $\theta_{1:n}$, governed by the modes of the mixing density $g$.  
The difference with respect to DP mixture models can be  explained %further appreciated 
by comparing the predictive rule $G_n$ as expressed by (\ref{eq:defGn}) with the predictive rule (\ref{eq:PolyaMixDP}) in DP mixtures. The comparison shows that the recursive rule $G_n$ implies a loss of information with respect to (\ref{eq:PolyaMixDP}), because it does not revise inference on $\theta_i$ as new data become available. The prior laws induced, respectively, by the recursive rule (\ref{eq:defGn})  and the \Polya sequence (\ref{eq:PolyaMixDP}), would be close if such a loss of information was negligible, that is, if $P(\theta_i \mid x_{1:i})$ was close to $P(\theta_i \mid x_{1:n})$. In general, this does not hold for small $i$ and $n$. 

\section{Asymptotic posterior laws} \label{sec:asymptoticLaw}

 By part (ii) of Theorem \ref{th:urnmart}, Newton's rule $G_n$ can be properly regarded as the point estimate, with respect to  quadratic loss, of the limit mixing distribution $G$ in a quasi-Bayes mixture model. Yet, it is desirable to go beyond point estimation, providing a full description of the uncertainty through the posterior distribution. % of $G \mid x_{1:n}$. 
We first obtain an asymptotic Gaussian approximation of the posterior distribution of $G(A)$, %given $x_{1:n}$, 
for any measurable set $A$. We then extend the results to the joint posterior distribution of a random vector $[G(A_1), \ldots, G(A_k)]'$.
%, for any measurable $A_1, \ldots A_k$. 

\subsection{Asymptotic posterior distribution and credible intervals.}

\label{sec:credint}
Let us recall that $P$ is a probability law for $(X_n)$ consistent with the assumptions (\ref{eq:newtonAsPred}). Here, we  present  an asymptotic Gaussian approximation of the conditional law $P(G(A) \in \cdot \mid x_{1:n})$, for a  measurable set $A$. 
For exchangeable sequences, central limit theorems and asymptotic results are usually given in terms of {\em stable convergence} (\cite{renyi1963}, \cite{aldous1985}, \cite{hausler2015}). The results below are in terms of {\em almost sure convergence of the conditional distributions}, or more briefly, {\em a.s. conditional convergence}, which is a stronger form of convergence (\cite{crimaldi2009}), that implies stable convergence and convergence in distribution of the unconditional law. 
Informally, Theorem \ref{th:central2} below, says that %, if $0<G(A)<1$, then 
\begin{equation} \label{eq:jointPosterior}
P( (G(A) - G_n(A)) \in \cdot \mid x_{1:n}) \approxeq \Norm(0, V_{A,n}\sum_{k>n} \alpha_k^2), 
\end{equation}
where 
%$\Norm(\mu,\sigma^2)$ denotes the Gaussian law with mean $\mu$ and variance $\sigma^2$, 
$V_{A,n}$ is defined in (\ref{eq:Vn1}), and the approximation holds for all $\omega=(x_1, x_2, \ldots)$ in a set of $P$-probability one. 
Notice that 
%the law $P$ describes an evolutionary process $(X_n)$ and 
 asymptotic results as the one above inform about the rate of convergence of the predictive probability $G_n(A)(x_{1:n})$ to the limit distribution $G(A)(x_1, x_2, \dots)$.
A novelty of our approach is in the statistical use we make of this kind of convergence, as informative of the asymptotic Gaussian form of the posterior distribution of the unknown $G(A)$. Indeed, we can read (\ref{eq:jointPosterior}) as 
$$ 
P( G(A) \in \cdot \mid x_{1:n}) \approxeq \Norm( G_n(A), V_{A,n}\sum_{k>n} \alpha_k^2), \quad \mbox{$P$-a.s.}
$$
Although having a similar flavor, these results differ from Bernstein-von-Mises types of theorems, which are stated a.s. with respect to the probability law  $P_{G_{true}}^\infty$ that assumes $X_i \iidsim F_{G_{true}}$. Our results are given a.s. with respect to the c.i.d. probability law $P$. Berstein-von-Mises results are a basic  tool for studying frequentist coverage of Bayesian procedures, which is beyond the scope of this paper. However, 
convergence of the kind (\ref{eq:jointPosterior}), too, provides insights on frequentist coverage properties, as we discuss in Section \ref{sec:simulation}.
%we provide some hints in the simulation study of Section \ref{sec:simulation}.  

%\medskip

We denote by $\Phi(t \mid \mu, \sigma^2)$ the distribution function of the $\Norm(\mu,\sigma^2)$ law, evaluated at $t$. A  $\Norm(0,0)$ is interpreted as the law degenerate at zero. 
Without loss of generality, we can assume that $f_{G_0}(x)\neq 0$ for every $x\in \mathbb X$. This implies 
\begin{equation} \label{eq:fGn-nonzero}
f_{G_n}(x)\neq 0\quad \mbox{for every $x\in\mathbb X$ and $n \geq 0$}. 
\end{equation}
Our first result finds a sequence  $(r_n)$ such that the conditional distribution of $\sqrt{r_n} (G(A)-G_n(A))$, given $X_{1:n}$, is asymptotically a zero-mean Gaussian law, with variance
\begin{equation}  \label{eq:V}
	   V_A = \int_{\{x:f_G(x) \neq 0\}} P_G(A\mid x)^2 
	dF_G(x)- G(A)^2,
\end{equation}
where, for any distribution function $H$ on $\Theta$, 
$P_{H}(A \mid x) =\int_A f(x\mid\theta) dH(\theta)/\int_\Theta f(x\mid\theta) dH(\theta)$. Before stating the theorem, we give the following lemma. 
Define, for any $A \in \mathcal B(\Theta)$ and $n \geq 1$, 
\begin{equation} \label{eq:Vn1}
V_{A,n} =
\int_{\mathbb X} P_{G_n}(A\mid x)^2
 dF_{G_n}(x) - G_n(A)^2. 
\end{equation}
Notice that $V_{A,n}$ can be written as $
V_{A,n}= E( (P_{G_n}(A \mid X_{n+1}) -  G_n(A))^2 \mid X_{1:n})$, expressing the prior-to-posterior variability, given $X_{1:n}$, when $G_n$ plays the role of the prior and $P_{G_n}(\cdot \mid x_{n+1})$ of the posterior.  

\begin{lemma} \label{lemma:V}
For any $A \in \mathcal B(\Theta)$, $ V_{A,n}$ converges to
$V_A$ $P$-a.s. as $n \rightarrow \infty$.   
\end{lemma}

We can now present the main results of this section. 

\begin{thm} 
	\label{th:rategeneral}
	Let $(\alpha_n)$ satisfy  $\sum_n\alpha_n=\infty$ and $\sum_n\alpha_n^2<\infty$ and let $(r_n)$ be a  monotone sequence of positive numbers  such that $r_n\sim (\sum_{k > n}\alpha_k^2)^{-1}$ as $n\rightarrow\infty$. If
\begin{equation}
\label{eq:centr_cond1}
	\sqrt{r_n}\sup_{k \geq n}\alpha_k\rightarrow 0
	\end{equation}
	 and
	 \begin{equation}
	 \label{eq:centr_cond2}
	  \sum_{k\geq 1}r_k^2\alpha_{k+1}^4<\infty, 
\end{equation}
	then, for every $A\in\mathcal B(\Theta)$, 
	\begin {equation}
	\label{eq:conv_gauss}
	P(	\sqrt{r_{n}} (G(A)-G_n(A)) \leq t \mid X_{1:n}) \rightarrow \Phi(t \mid 0,  {V_{A}}) \quad\mbox{ $P$-a.s.,}
	\end{equation}
	with ${V_{A}}$ as in (\ref{eq:V}).
	If $\alpha_n=(\alpha+n)^{-\beta}$ with $1/2 <\beta\leq 1$ and $\alpha>0$, then 
	%If 	$\alpha_n=1/(\alpha+n)$ with  $\alpha>0$, then, 
	(\ref{eq:conv_gauss}) holds with 
	$r_n=(2\beta-1)n^{2\beta-1}$.
	\end{thm}

\begin{rmrk}
Assumptions (\ref{eq:centr_cond1}) and (\ref{eq:centr_cond2}) hold for most choices of $(\alpha_n)$ satisfying $\sum_n\alpha_n=\infty$ and $\sum_n\alpha_n^2<\infty$. In particular, if  $(\alpha_n)$ is ultimately decreasing, then (\ref{eq:centr_cond1}) is a consequence of (\ref{eq:centr_cond2}). A sufficient condition for (\ref{eq:centr_cond2}) is $ \alpha_n=1/(n b_n)$ for a sequence $b_n$ which is ultimately non increasing. Indeed, in this case 
$$
	\limsup_{n\rightarrow\infty}r_{n-1}\alpha_n
	=\limsup_{n\rightarrow\infty} \frac{(nb_n)^{-1}}{\sum_{k\geq n}(kb_k)^{-2}}\leq 
	\limsup_{n\rightarrow\infty} \frac{(nb_n)^{-1}}{b_n^{-2}\sum_{k\geq n}k^{-2}}\leq\limsup_{n\rightarrow\infty}b_n<\infty.
		$$ 
	In turn, this implies that  $r_n^2 \alpha_{n+1}^4<(\sup_kb_k)^2\alpha_{n+1}^2$, for $n$ large enough, and therefore (\ref{eq:centr_cond2}).
	\end{rmrk} 

\begin{rmrk} If $\omega=(x_1, x_2, \ldots)$ is such that $V_A(\omega)=0$, then Theorem \ref{th:rategeneral} ensures convergence to a degenerate distribution on zero.
From the definition of $V_A$, 
 it is immediate to see that $V_A(\omega)=0$ if and only if $P_G(A \mid X)(\omega) = G(A)(\omega)$, which happens if and only if $G(A)(\omega)$ is zero or one. 
\end{rmrk} 

In Theorem \ref{th:rategeneral}, the limit variance $V_A$ depends on  $G$ and is, therefore,  unknown. By Lemma \ref{lemma:V}, a convergent estimator is given by $V_{A,n}$. In i.i.d. settings, one would exploit Cram\'{e}r-Slutzky Theorem to replace the random $V_A$ with its consistent estimate, obtaining an asymptotic distribution that allows to compute asymptotic credible intervals for $G(A)$. In the present case, however, this is not immediate, as $V_A$ is random; moreover, we want to prove convergence of the {\em conditional} distributions. 
For the unconditional distribution, the presence of a random quantity in the limit is solved through stable convergence. Dealing with the conditional distributions, we need further work; we first prove convergence of the joint conditional distribution of $(\sqrt{r_n} (G_n(A)-G(A)), V_{A,n})$ given $x_{1:n}$ to finally prove the following Theorem.

\begin{thm} 	\label{th:central2}
Let $A \in \mathcal B(\Theta)$. Then, with $(r_n)$ defined as in Theorem \ref{th:rategeneral} and under the same assumptions, 
for $P$-almost all $\omega=(x_1,x_2,\dots)$ such that $V_A(\omega)>0$,
\begin{equation}
\label{eq:conv_2}
P(\sqrt{r_{n}} \frac{G(A)-G_n(A)}{\sqrt{V_{A,n}}}\leq t \mid  x_{1:n}) \rightarrow \Phi(t \mid  0,1). \end{equation}
If $\alpha_n=(\alpha+n)^{-\beta}$ with $1/2 <\beta\leq 1$ and $\alpha>0$, then 
(\ref{eq:conv_2}) holds with 
$r_n=(2\beta-1)n^{2\beta-1}$.
\end{thm}

%\medskip

 Theorems \ref{th:rategeneral} and \ref{th:central2} allow to obtain  asymptotic credible intervals for G(A). Indeed, for a fixed set $A$, it follows from
 Theorem \ref{th:central2} that, for $P$-almost all $\omega=(x_1, x_2, \ldots)$ such that $V_A(\omega)>0$,
$$
P( G_n(A) - z_{1-\gamma/2} \sqrt{V_{A,n}/r_n} < G(A) < G_n(A) +  z_{1-\gamma/2} \sqrt{ V_{A,n}/r_n} \mid x_{1:n}) \approx 1-\gamma,  
$$
where $z_{1-\gamma/2}$ is the $(1-\gamma/2)$-quantile of the standard Gaussian distribution. If $V_A(\omega)=0$,
then Theorem \ref{th:rategeneral} implies that the limit distribution is degenerate on zero, therefore, for any $\epsilon > 0$
$$
P( G_n(A) - z_{1-\gamma/2} \sqrt{ \epsilon/r_n} < G(A) < G_n(A) +  z_{1-\gamma/2} \sqrt{\epsilon/ r_n} \mid x_{1:n}) \geq 1-\gamma, 
$$
asymptotically. It follows that, for every $\epsilon >0$, 
$$
\left[ G_n(A) - z_{1-\gamma/2} \sqrt{\frac{max ( V_{A,n}, \epsilon)}{r_n}} ; G_n(A) +  z_{1-\gamma/2} \sqrt{\frac{max ( V_{A,n}, \epsilon)}{r_n}} \right]
$$ is an asymptotic credible interval for $G(A)$, of level at least $1-\gamma$.

\subsection{Asymptotic joint posterior distribution and credible regions.}
\label{sec:credreg}
 
We now study the joint behavior of $(G_n(A_1)-G(A_1),\dots,G_n(A_k)-G(A_k))$, for any fixed choice of $A_1,\dots,A_k \in\mathcal B(\Theta)$.  
 As in the previous section, we assume that $f_{G_0}(x)\neq 0$ for every $x\in \mathbb X$, which implies (\ref{eq:fGn-nonzero}).   
For every $n\geq 1$, and $A_i,A_{i'}$, let 
\begin{eqnarray*}
C_{A_i,A_{i'},n}&=Cov\left( P(\theta_{n+1}\in A_i\mid 
X_{1:n+1}),P(\theta_{n+1}\in A_{i'}\mid  X_{1:n+1})\mid X_{1:n}\right)\\
&=\int_{\mathbb X}P_{G_n}(A_i \mid x)P_{G_n}(A_{i'} \mid x)dF_{G_n}(x)
-G_n(A_i)G_n(A_{i'}), 
\end{eqnarray*}  
and 
$$
C_{A_i,A_{i'}}=\int_{\{x:f_G(x)\neq 0\}}P_{G}(A_i \mid x)P_{G}(A_{i'} \mid x)dF_{G}(x)-G(A_i)G(A_{i'}).
$$
Following the same line of reasoning as in Lemma \ref{lemma:V}, it can be proved that, as $n\rightarrow\infty$, 
 $$C_{A_i,A_{i'},n}\rightarrow C_{A_i,A_{i'}} \quad \mbox{$P$-a.s.,}$$
and that, denoting by $C_n(A_1,\dots,A_k)$
the matrix $\left[C_{A_i,A_{i'},n}\right]_{i,i'}$,
\begin{equation} \label{eq:VV2}
C_n(A_1,\dots,A_k) \rightarrow C(A_1,\dots,A_k)\equiv \left[C_{A_i,A_{i'}}\right]_{i,i'} \quad \mbox{$P$-a.s.}
\end{equation}
We denote by $\Phi_p(\mathbf{t} \mid \bm{\mu}, \Sigma)$ the distribution function of the $p$-dimensional Gaussian law $\Norm_p({\bm\mu},\Sigma)$ with mean vector $\bm{\mu}$ and covariance matrix $\Sigma$.
Then, we can prove the following theorem.
\begin{thm}
	\label{th:ratevett1}
In model (\ref{eq:newtonAsPred}), let $(\alpha_n)$ satisfy  $\sum_n\alpha_n=\infty$ and $\sum_n\alpha_n^2<\infty$.
Let $(r_n)$ be a monotone sequence of positive numbers such that $r_n\sim(\sum_{j>n}\alpha_j^2)^{-1}$ as $n\rightarrow\infty$.
If (\ref{eq:centr_cond1}) and (\ref{eq:centr_cond2}) hold, then,
for every $k \geq 1$ and every $A_1,\dots,A_k\in\mathcal B(\Theta)$, 
\begin{equation}
\label{eq:conv_joint}
P (\sqrt{r_n} 
\begin{bmatrix}
G(A_1)-G_n(A_1) \\
\vdots\\
G(A_k)-G_n(A_k) 
\end{bmatrix} 
\in \cdot \mid X_{1:n} ) \rightarrow \Phi_k( \cdot \mid {\bf 0},  C(A_1,\dots,A_k)), \quad\mbox{$P$-a.s.}\mbox{ for $n \rightarrow \infty$},
\end{equation}
with $C(A_1,\dots,A_k)$ as in (\ref{eq:VV2}).
	If $\alpha_n=(\alpha+n)^{-\beta}$ with $1/2 <\beta\leq 1$ and $\alpha>0$, then 
	(\ref{eq:conv_joint}) holds with 
	$r_n=(2\beta-1)n^{2\beta-1}$. 
\end{thm}

The following result is the analogous of Theorem \ref{th:central2} for the joint posterior distribution.
\begin{thm}
	\label{th:ratevett2}
	Under the same assumptions as in Theorem \ref{th:ratevett1}, for every $k \geq 1$ and every $A_1,\dots,A_k\in\mathcal B(\Theta)$, 
\[
P\left(\sqrt {r_n}\; C_n(A_1,\dots,A_k)^{-1/2}
\begin{bmatrix}
G(A_1)-G_n(A_1)\\
\dots\\
G(A_k)-G_n(A_k)
\end{bmatrix}
\in \cdot  \mid x_{1:n} \right)\rightarrow \Phi_k(\cdot \mid {\bf 0}, I) \quad \mbox{for $n \rightarrow \infty$},
\]
for $P$-almost all $\omega=(x_1, x_2, \ldots)$ such that $
\det(C(A_1,\dots,A_k)(\omega)) \neq 0$. 
\end{thm}

Based on Theorems \ref{th:ratevett1} and \ref{th:ratevett2}, we can find an asymptotic credible region for $(G(A_1),\dots,G(A_k))$.

\begin{prop} \label{th:crebibleRegion}
Let  $\mathbf G_n(\mathbf A)=[G_n(A_1), \dots, G_n(A_k)]'$ and let $\chi^2_{1-\gamma}$ denote the $(1-\gamma)$-quantile of the chi-square distribution with $k$ degrees of freedom. Then, for every $\epsilon >0$, the set
\[
E_n^{(\epsilon)}=\left\{\mathbf s\in\mathbb R^k:(\mathbf s-\mathbf G_n(\mathbf A))'
( C_n(A_1,\dots, A_k)+\epsilon I)^{-1}(\mathbf s-\mathbf G_n(\mathbf A))
\leq \frac{\chi^2_{1-\gamma}}{r_n}\right\}
\]
satisfies, $P$-a.s., 
\[
\liminf_n P([G(A_1),\dots,G(A_k)]'\in E_n^{(\epsilon)}\mid
X_{1:n})\geq 1-\gamma.
\]
\end{prop}

\section{Recursive prediction and learning} \label{sec:estensioni}

Once the statistical modeling assumptions underlying the recursive rule $g_n$ are clear,  the role of the different ingredients of the model becomes clearer, too, and several statistical applications can be  naturally envisaged. 
We first describe a c.i.d. time-varying mixture model that is consistent with the assumptions (\ref{eq:newtonAsPred}) and has a natural statistical interpretation,  underlining the temporal nature of the data implied by the recursions. 
This model may be of independent interest as a model for temporal data. 
%In this case, the results of the previous section provide an exact Bayesian method for recursive learning and prediction. 
%estimating the limit mixing density, with proper uncertainty quantification.
When used, instead, as a fast approximation of a (static) exchangeable mixture model, the time-dependent mixture  specification gives further intuition on the role of the parameters, in particular of the weights $\alpha_n$, that control the dynamics of the model. 

\subsection{A time-dependent mixture model} \label{sec:timeMixture}
%In a standard mixture model, as previously discussed, $X_i \mid G \iidsim f_G=\int f(x \mid \theta) dG(\theta)$. 
In the quasi-Bayes mixture model (\ref{eq:newtonAsPred}), the sequence $(X_n)$ is not exchangeable and, therefore, there is no random distribution $\tG$, such that $X_i \mid \tG \iidsim f_{\tG}$.  We can, however, think of a sequence of latent random distributions, say $(\tG_n)$, such that 
\begin{equation} \label{eq:modLatentGn}
X_n \mid \tG_n \indsim F_{\tG_n}(\cdot)=\int F(x \mid \theta) d\tG_n(\theta).
 \end{equation} 
This is a time-varying mixture model that envisages a temporal evolution of the latent distributions $\tG_n$. For example, the state $\theta_n$ could express student's skill and $\tG_n$ the distribution of skill in the class at time $n$.  
Or one may want to model an imbalance due to some intervention in the system under study (e.g., in ecologic studies, an imbalance in the population of species due to climate or to some form of competition), which breaks the symmetry of an exchangeable setting. 

Models of this kind have been considered in the Bayesian literature for nonparametric density estimation with temporal data, usually assigning dependent DP priors on the $\tG_n$.  Here, we specify a type of unpredictable dynamics of the latent sequence  $(\tG_n)$, 
made precise in the following proposition, such that $E(\tG_n)=G_0$ for any $n$ and the $\tG_n$ converge to a random limit distribution $\tG$. Informally, one models an imbalance in the system, that tends to stabilize in the long run, so that the system converges to a (new) stationary steady state. The resulting process is c.i.d.

\begin{prop} \label{prop:statespace}  
	Suppose that $Y_n\mid\tilde H_n\stackrel{indep}{\sim} \tilde H_n$. If, for every $n \geq 1$, $E(\tilde{H}_n)=H_0$  and 
\[
E(\tilde H_{n+2}\mid y_{1:n})=E(\tilde H_{n+1}\mid y_{1:n}),\] 
then   $(Y_n)$ is c.i.d. with directing random measure 
$\tilde H=\lim_n E(\tilde H_n\mid y_{1:n-1})$, $P$-a.s. \\
In particular, if $\tilde H_n-E(\tilde H_n\mid y_{1:n-1})\rightarrow 0$, $P$-a.s., then
	$\tilde H=\lim_n \tilde H_n$, $P$-a.s.
\end{prop}  

In this setting, let us assume the following dynamics for the latent sequence  $(\tG_n)$ in the temporal mixture model (\ref{eq:modLatentGn}) 
\begin{eqnarray}  \label{eq:tG2}
  \tG_1 &\sim &  \DP((1-\alpha_1)/\alpha_1 \;G_0) \nonumber \\
\tG_n  \mid X_{1:n-1}, \theta_{1:n-1} &\sim & \DP((1-\alpha_n)/\alpha_n\; G_{n-1}),  \quad \quad n > 1.
\end{eqnarray}
That is, the conditional law of $\tG_n$ is a Dirichlet process centered on the current estimate $G_{n-1}$, computed according to Newton's one-step-ahead updating rule (\ref{eq:newton}). 
This dynamics is consistent with the modeling assumptions (\ref{eq:newtonAsPred}). 
Notice that, as $x_n$ becomes available, the updating (\ref{eq:newton}) is exact, that is, it is indeed  the {\em Bayesian} point estimate  $E(\tG_n \mid  x_{1:n})$ of $\tG_n$ from the DP prior (\ref{eq:tG2}). 
It is easy to verify that the sequence $(F_{\tG_n})$, with $(\tG_n)$ as in (\ref{eq:tG2}), satisfies the assumptions of Proposition \ref{prop:statespace} and that $F_{\tilde G_n}$ converges $P$-a.s. to $F_{{G}}$, with $G=\lim_n G_n$. Thus, by Proposition \ref{prop:hconv}, the sequence $(X_n)$ is c.i.d. and, asymptotically, $X_n \mid G \approxiid F_{G}$.

\begin{figure}[t!]  
 \centering 
  \includegraphics[width=.8\textwidth, height=8cm]{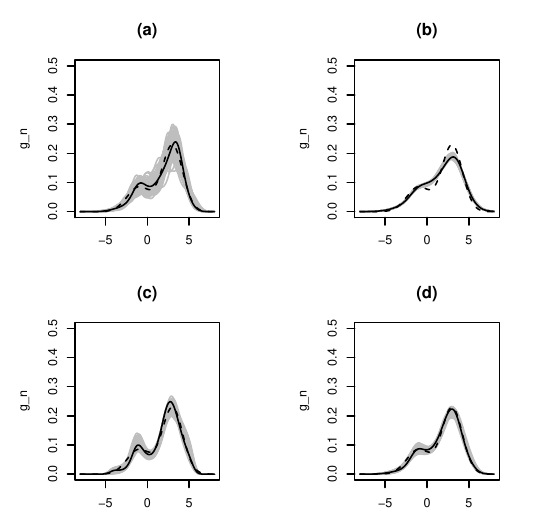}%width=\linewidth]
\caption{Mixing density estimate $g_n$, and estimates obtained over $100$ random permutations of the original sample (plotted in gray). Simulated data from a location mixture of Gaussians; $\sigma^2=1; n=1000$. The true mixing density is the dashed curve. % $0.3 \Norm(-1, 2) + 0.7 \Norm(3, 1.5)$. 
%Initial distribution  $G_0=\Norm(1,9)$. 
Panel (a): $\alpha_n=1/(\alpha+n)$, with $\alpha=1$. 
Panel (b): $\alpha_n=1/(\alpha+n)$, $\alpha=100$. 
Panel (c): $\alpha_n=1/(\alpha+n)^{2/3}$, $\alpha=100$. 
Panel (d): split-sample weights, $N=500$, $\gamma=3/4$; $\alpha=100$.}
\label{fig:bimodal} 
\end{figure}

\subsection{Recursive learning on the asymptotic mixing density} 
\label{sec:simulation}
When  the time-varying mixture model (\ref{eq:modLatentGn})-(\ref{eq:tG2}) is used as the actual model for temporal data,  the results in the previous sections give a fully Bayesian method for recursive learning and prediction. 
%Bayesian online prediction of $\theta_{n+1}$ given $x_{1:n}$ is solved through the predictive density $g_n$ of $\theta_{n+1} \mid x_{1:n}$, that can be computed recursively according to  (\ref{eq:recursiveDensities}). Bayesian recursive learning on the limit  $G(t)$ of the $\tG_n(t)$ is solved through the posterior distribution of $G(t)$, and  Theorems \ref{th:rategeneral} and \ref{th:central2} provide an asymptotic  Gaussian approximation of it. 
%%%%
%%%% INIZIO FIGURE MLTINODALE E UNIMODALE
%%%%%%%%%%%%%%%%%%%%%%%%%%%%%%%%%%%%%%%%%%%%%%%%%
%%%%%%%%%%%%        MULTIMODALE       %%%%%%%%%%%%%%%%%%%%%
%%%%% Dati salvati come xEx1multimodale (per n=5000) e xEx1Multimodale.n10mila %%%%%%
\begin{figure}[ht]  
\centering    
\includegraphics[width=.8\textwidth, height=8cm]{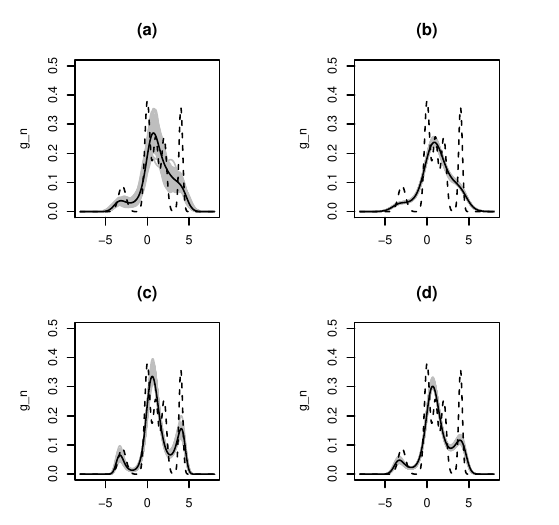} 
\caption{Mixing density estimate $g_n$ (black) and estimates obtained over $2000$ random permutations of the original sample (plotted in gray). 
Simulated data from a location mixture of Gaussians: $\sigma^2=1, n=5000$ and multimodal mixing density (dashed curve).  
%$g^*= \sum_{j=1^5} p_j \Norm(\mu^*_j, \tau_j^*)$ with $(p_1, \ldots, p_5)=(.1, .3, .2, .2, .2), (\mu^*_1, \ldots, \mu^*_5)= (-3, 0, 2, 1, 4), (\tau^*_1, \ldots, \tau_5^*)=(.2, .1,.1, .1, .05)$. 
%Initial $G_0$ and 
Weights $\alpha_n$ as in Figure \ref{fig:bimodal}. } 
\label{fig:multimodale-s1}
\end{figure} 
%%%%%%%%%%%%%%%%%%%%%%%%%%%%%%
%  Stesssa cosa, ma con $\sigma^2= 0.1$: 
%%% DATI SALVATI : xEx2Multimodale
\begin{figure}[ht]   
\centering
\includegraphics[width=.8\textwidth, height=8cm]{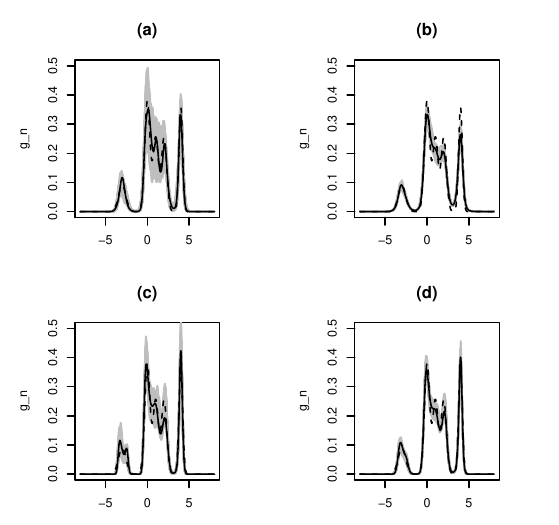} 
\caption{Mixing density estimate $g_n$ (black) and estimates obtained over $200$ random permutations of the sample (gray). Simulated data from a location mixture of Gaussians: $\sigma^2=0.1, n=5000$. Multimodal mixing density (dashed) and weights $\alpha_n$ as in Figure \ref{fig:bimodal}.} 
%$g^*= \sum_{j=1^5} p_j \Norm(\mu^*_j, \tau_j^*)$ with $(p_1, \ldots, p_5)=(.1, .3, .2, .2, .2), (\mu^*_1, \ldots, \mu^*_5)= (-3, 0, 2, 1, 4), (\tau^*_1, \ldots, \tau_5^*)=(.2, .1,.1, .1, .05)$. Initial $G_0$ and weights $\alpha_n$ as in Figure \ref{fig:multimodale}.}
%\end{center}
\label{fig:multimodale-s01}
\end{figure} 
%%%%%%%%%%%%%%%%%%%%%%%%%%%
%% UNIMODALE, sigma^2=5   -- %%FIGURA OK, TAGLIATA
%\begin{figure}[ht]   
%\centering
%\includegraphics[width=0.8\textwidth]{Rplot-unimodale-s5-alpha50.pdf} 
%\caption{Mixing density estimate $g_n$ (black) and estimates $g_n^{(\pi)}$ obtained over $50$ random permutations $\pi$ of the original sample (plotted in gray). Simulated data from a location mixture of Gaussians: $\sigma^2=5, n=1000$, unimodal true mixing distribution (dashed curve) $\Norm(2,2)$. Initial distribution $G_0=\Norm(1,9)$. Panel (a): $\alpha_n=1/(\alpha+n)$ with $\alpha=5$. Panel (b): $\alpha_n=1/(\alpha+n)$,  $\alpha=50$. Panel (c): $\alpha_n=1/(\alpha+n)^{2/3}$, $\alpha=50$. Panel (d): Split-sample weights, $N=500$, $\gamma=3/4$; $\alpha=50$.} 
%\label{fig:unimodale-s5-alpha50.pdf}
%\end{figure} 
%%%%%%%%%%%
%%%%%%%%% 
%% UNIMODALE, sigma^2=.1, due diverse $g^*$
%% Per g^*=N(2,2): dati salvati come xEx1Unimodal.N22
%% Per g^*=N(2,0.2):  dati salvati come xEx1Unimodal.N202)
\begin{figure}[ht]   
\centering
\includegraphics[width=0.6\textwidth,height=5cm]{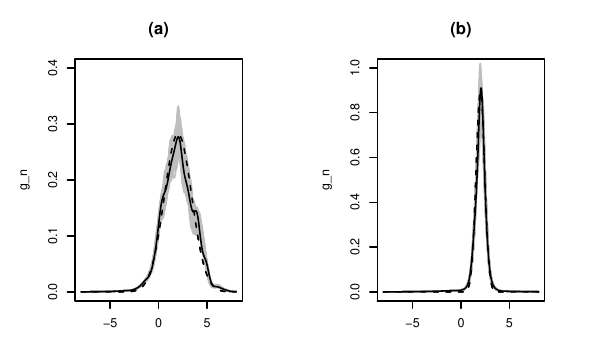}   
%{Rplot-unimodale-s5-alpha50.pdf} 
\caption{Mixing density estimate $g_n$ (black) and estimates for $100$ random permutations of the original sample (gray). Simulated data from a location mixture of Gaussians: $\sigma^2=0.1, n=1000$.  %Initial distribution $G_0=\Norm(1,9)$. 
Panel (a) Mixing density (dashed) $g^*=\Norm(2,2)$. 
Panel (b)  Mixing density (dashed) $g^*=\Norm(2,0.2)$.
Weights $\alpha_n = 1/(\alpha+n)$, with $\alpha=50$. 
}
%\end{center}
\label{fig:unimodale-differentgtrue.pdf}
\end{figure} 
%%%%%%%%%%%%%%%%%%%%%%%%%  FINE FIGURE %%%%%%%%%%%%%%%
% 
In a static setting, (as in previous sections), it may also be regarded as a misspecified but computationally fast model that approximates an exchangeable mixture model.
 While with temporal data the dependence of the prediction of $\theta_{n+1} \mid x_{1:n}$ on the ordering of the observations is natural, in this case sensitivity to the ordering is a drawback of the lack of exchangeability, which one wants to attenuate. With this aim, we discuss the role of the various ingredients of the model via a simulation study.  In particular, we notice that the weights $\alpha_n$ have a dual role,  controlling the speed of convergence to exchangeability as well as the learning rate of the predictive rule. Roughly speaking, one has approximate  exchangeability when the predictive distribution $G_n$ is close to its limit $G$ (\cite{aldous1985}, Lemma 8.2; see  Section \ref{sec:quasi-Bayes}). The results in Section \ref{sec:asymptoticLaw} provide the rate of convergence of $G_n$. In the class of weights of the form $\alpha_n=1/(\alpha+n)^\beta$, with $\beta \in (1/2, 1]$, the fastest convergence rate $1/\sqrt r_n=1/\sqrt n$ is obtained for $\beta=1$. The time-dependent mixture-model provides further intuition. Here, the $\alpha_n$ affect the dynamics of the random distributions $\tG_n$. Roughly speaking, smaller $\alpha_n$ give a milder evolution of the $\tG_n$, and thus a situation closer to exchangeability, for which the $\tG_n$ would remain constant.   

However, there is a subtle trade off: on one hand, weights $\alpha_n$ that rapidly decay to zero allow to quickly reach asymptotic exchangeability; on the other hand, the  $\alpha_n$ determine the weight of the current observation in the predictive distribution $G_n$, and small values may lead to poor learning. Thus, one may want to use small weights $\alpha_n$  that  do not however decrease to zero too quickly. This could be obtained by using weights of the form $\alpha_n=1/(\alpha+n)^{\beta_n}$, thus letting the exponent depend on $n$. 
A practical suggestion is to split the sample $x_{1:n}$, using $\beta_n=1$, with a fairly large $\alpha$, for an initial {\em prior-training} sample, say $n \leq N$, in order to rapidly reach approximate exchangeability. Then use weights that decrease to zero  slowly, thus with $\beta_n=\gamma < 1$, for the {\em learning sample}, i.e. for  $ n>N$, in order to more efficiently learn from the data, once in a situation of approximate exchangeability. For brevity, let us denote this choice as {\em split-sample weights} with parameters $N, \gamma$.  
 
The lack of exchangeability of the recursive estimate is well known, and is addressed in the literature by using an average of the estimates obtained over a number of random permutations of the original sample; see for example \cite{tokdarMartinGhosh2009}. For a fixed sample size, this procedure is still very fast.  But, when observations arrive sequentially, the recursive feature of the computations  is lost and the complete procedure must be re-initiated each time a new observation becomes available. The computational cost may still be reasonable for a fairly small number $M$ of permutations, that already gives improved results. However, the results depend on $M$.  
Moreover, the total number of permutations of $x_{1:n}$ rapidly increases with $n$. If one wants the proportion $M/n!$ of visited permutations to be constant, %the choice of 
then $M$ should  increase with $n$. The recursive estimate (\ref{eq:newton}) remains computationally attractive. Understanding the role of its components  is important in this permutation-based setting, as a smaller value of $M$ is needed when sensitivity to permutations is attenuated by an adequate choice of $\alpha_n$. 
  
In the following examples,  the data are generated from a location mixture of Gaussian distributions; that is, we generate $\theta_i$, $i=1, \ldots, n$, i.i.d. from a mixing density $g^*$, and $X_i\nobreak\mid\nobreak\theta_i\nobreak\indsim\nobreak\Norm\nobreak(\theta_i,\nobreak\sigma^2\nobreak)$, $i=1, \ldots, n$, with $\sigma^2$ known. 
We start from a vague initial distribution $G_0=\Norm(1,9)$. 
In this case, the assumptions of Theorem \ref{th:absolcont} hold, therefore the prior law under the c.i.d model  selects $P$-a.s. absolutely continuous distributions. 
As our results are asymptotic, we consider fairly large values of $n$.

We run the simulation for different choices of the weights $\alpha_n$ and different shapes of the  mixing density $g^*$. The results 
all lead to similar conclusions. We show the results for a bimodal, a multimodal and a unimodal true mixing density $g^*$, and compare weights of the form $\alpha_n = 1/(\alpha+n)^\beta$, with $\beta \in (1/2, 1]$, for different choices of $\alpha$ and $\beta$. 
  
In Figure \ref{fig:bimodal}, the true mixing density $g^*$ (dashed curve) is bimodal, a mixture of two Gaussian densities $g^*= 0.3 \Norm(-1,2)+0.7 \Norm(3, 1.5)$.  The sample size is $n=1000$. 
We plot the %Sub-figure (\ref{fig:bimodal}a) 
%\ref{fig:bimodal-gn} shows the 
recursive mixing density estimate $g_n$, for different choices of the weights $\alpha_n$, together with the estimates (plotted in gray) 
%In sub-figure (\ref{fig:bimodal}b), %\ref{fic:bimodal-perm}, 
% the black curve is the average of the estimates (plotted in gray) 
obtained over $100$ random permutations of the original sample $x_{1:n}$ (the same $x_{1:n}$ and the same permutations in each panel), to give %for giving 
an idea of sensitivity to permutations. 
In panel (a) 
%of both sub-figures,  
we use the popular choice of DP-like weights $\alpha_n=1/(\alpha+n)^\beta$ with $\beta=1$ and a small value of $\alpha=1$. Sensitivity of the estimates to the ordering of the observations is evident. The reason is that, for $\alpha=1$, the weight $\alpha_n$ is too big when $n$ is small. 
In panel (b), we consider a large value of $\alpha=100$; this choice gives a small weight $\alpha_n=1/(\alpha+n)$ also for small values of $n$. The effect of the ordering is greatly attenuated, but the model does not learn enough from the data, suggesting that the weights are too rapidly decreasing to zero. 
Panel (c) shows the estimate $g_n$ obtained using $\beta=2/3$ as often suggested in the literature (\cite{martinTokdar2009}, \cite{dixitMartin2019}), with $\alpha=100$. These weights appear to give  good learning, yet they decay to zero too slowly, which again affects sensitivity to the ordering.  
Finally, in panel (d), we let $\beta$ depend on $n$, using the simple split-sample weights, with $N=500$, $\gamma=3/4$, and $\alpha=100$. The c.i.d. model 
%Model (\ref{eq:modLatentGn})-(\ref{eq:tG2}) 
remains misspecified (not  exchangeable), but now the effect of the ordering is reduced and the learning rate is fairly satisfactory. Again, the advantage of the (slightly) misspecified model is the speed of the  recursive computations. 

We repeat the simulation with a multimodal mixing density $g^*= \sum_{j=1^5} p_j \Norm(\mu^*_j, \tau_j^*)$ with $(p_1, \ldots, p_5)=(.1, .3, .2, .2, .2)$;  $(\mu^*_1, \ldots, \mu^*_5)= (-3, 0, 2, 1, 4)$ and $(\tau^*_1, \ldots, \tau_5^*)=(.2, .1,.1, .1, .05)$. Again we fit a location mixture of Gaussian distributions, with $\sigma^2=1$. 
The sample size is $n=5000$. The recursive density estimate $g_n$ is shown in Figure \ref{fig:multimodale-s1}, together with the estimates (plotted in gray) obtained over $200$ random  permutations
of the original sample. The true mixing density $g^*$ is the dashed curve. The weights $\alpha_n$ are as in Figure \ref{fig:bimodal}. In all cases, the estimate $g_n$ reasonably reconstructs the bulk of the masses of $g^*$, but 
oversmooths the three central modes. In Figure \ref{fig:multimodale-s01}, the simulation setting is as before, but now $\sigma^2=0.1$. Not surprisingly, better results are obtained for a smaller kernel variance $\sigma^2$.  The estimates are improved, with the choice of split-sample weights in panel (d) proving a good compromise between attenuated sensitivity to permutations and efficient learning.  

As illustrated in the previous example, the kernel variance $\sigma^2$ affects the smoothness of the estimate $g_n$. Indeed, a small value of $\sigma^2$ tends to give a component $P(\theta_i \mid x_{1:i})$ in the predictive distribution (\ref{eq:defGn}) which concentrates  around $x_i$, thus favoring picks which track the $x_i$, and accentuating the sensitivity of $g_n$ to the ordering of the data. 
However, there is trade off with the dispersion of the $\theta_i$ generated from $g^*$. If $g^*$ is quite concentrated and $\sigma^2$ is small, the $x_i$ tend to be concentrated, too, thus attenuating  the data-tracking, order-dependent behavior. 
In Figure  \ref{fig:unimodale-differentgtrue.pdf}, the data are simulated from a location mixture of Gaussian kernels with a fairly small variance $\sigma^2=0.1$; the sample size is $n=1000$.  
For brevity, we only report the recursive estimates  $g_n$  obtained for DP-like weights $\alpha_n=\alpha/(\alpha+n)$, with $\alpha=50$ (black curve). 
In panel (a), the true mixing density is a $\Norm(2,2)$. The estimate $g_n$ gives a  reasonable idea of the shape of $g^*$, but it is too wiggly. A large value of $\alpha$ may be used to obtain a smoother estimate. In panel (b), the true mixing density $g^*$ is $\Norm(2, 0.2)$. As expected, with a more concentrated $g^*$, the estimate $g_n$ is smoother, and less sensitive to permutations.

\subsection{Inference on $G$}  \label{sec:inferenceG}

%%% intervalli credibilita' 
% Bimodale. DATI salvati: xEx1Normale
\begin{figure}[t!] 
%\begin{center}
\centering
\includegraphics[width=0.8\textwidth, height=7cm]{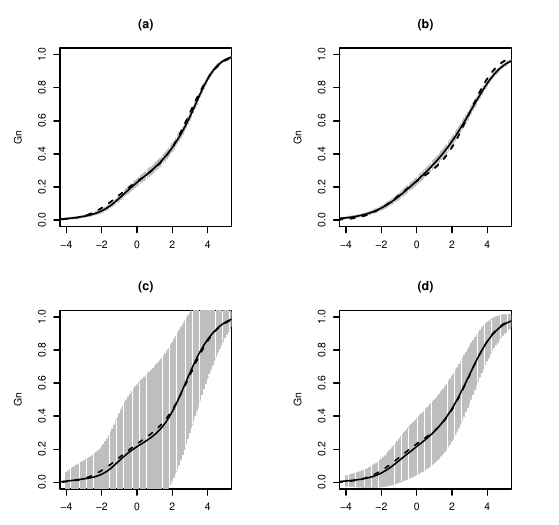}
%{Rplot-s1-alpha100-G-NEW2.pdf}  
%{plot-Gn-new2} %{Rplot42.pdf}
\caption{Recursive estimate $G_n$ (solid curve) and asymptotic $95\%$  marginal credible intervals. Dashed curve: true mixing distribution. 
%Location mixture of Gaussians; $\sigma^2=1; n=1000$. 
Data, weights $\alpha_n$ and bimodal mixing density as in Figure \ref{fig:bimodal}. 
 % as evidenced in the plot. 
%\red{Dotted curve: mean over $50$ random permutations}. 
}
%\ref{fig:order-bimodal}.}   
\label{fig:Gn1000}
\end{figure}
%
% Multimodale. DATI: xEx1Multimodale
\begin{figure}[th]
\centering
\includegraphics[width=0.8\textwidth, height=7cm]{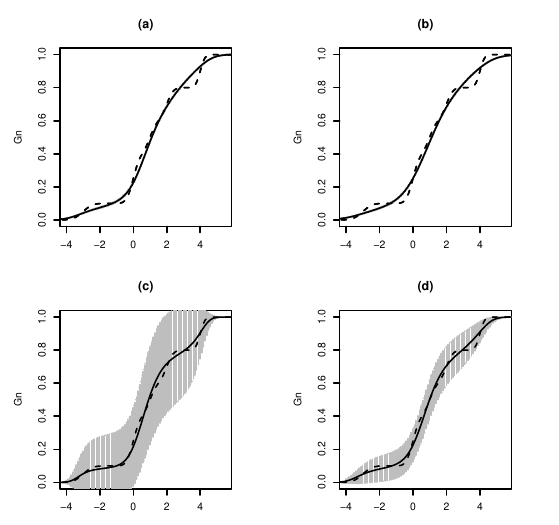}  
\caption{ Recursive estimate $G_n$ (solid curve) and asymptotic $95\%$  marginal credible intervals. % for $G(t_j)$, on a fine grid of values $t_j$. 
Dashed curve: true mixing distribution.
%Location mixture of Gaussians; $\sigma^2=1; n=5000$. 
Data, weights $\alpha_n$ and multimodal mixing density as in Figure \ref{fig:multimodale-s1}.}
\label{fig:GnMultimodale}
\end{figure}

Let us now consider inference on the mixing distribution. For any $t$, the recursive rule (\ref{eq:newton}) provides a point estimate $G_n(t)=E(G(t) \mid x_{1:n})$, for any $n \geq 0$. Moreover, we can use the asymptotic approximation of the posterior distribution of $G(t)$, given in Section \ref{sec:asymptoticLaw}, to provide asymptotic marginal  credible intervals. 
Figure \ref{fig:Gn1000} shows the results for the same data $x_{1:n}$ and the same choice of the weights $\alpha_n$ as in Figure \ref{fig:bimodal}. The true mixing distribution  $G^*$ (dashed curve) corresponds to the bimodal density $g^*$ in Figure \ref{fig:bimodal}. 
%\red{ Similar results are obtained for the case of a unimodal mixing density as in Figure \ref{fig:unimodale-s5-alpha50.pdf}.}
We show the point estimate $G_n(t)$ (solid curve), together with the $95\%$ asymptotic marginal credible intervals (gray), for $t$ in a fine grid. %grid as evidenced in the plot. 
Notice that, as the simulated data are i.i.d. from a known mixture density $f_{G^*}$, where $G^*$ plays the role of the true mixing distribution, this study provides some intuition on frequentist coverage of the quasi-Bayes procedure.  

The expression of the asymptotic intervals was obtained in Section \ref{sec:credint}. For weights $\alpha_n=1/(\alpha+n)^\beta$ with $\beta \in (1/2, 1]$, 
$$ G(t) \mid x_{1:n} \approx \Norm (G_n(t), \frac{V_{(-\infty, t], n}}{r_n}) $$
where $r_n=(2 \beta - 1) n^{2 \beta-1}$. Thus, the estimated asymptotic variance is 
$V_{(-\infty, t], n}/n$ for $\beta=1$ (panels (a) and (b));  $ 3 V_{(-\infty, t], n}/n^{1/3}$ for $\beta=2/3$ (panel (c)) and $ 2 V_{(-\infty, t], n}/n^{1/2}$ for $\beta=3/4$ (panel (d)). 

The results complement the discussion on Figure \ref{fig:bimodal}.
A choice of the weights $\alpha_n=1/(\alpha+n)^\beta$ with $\beta=1$ gives the fastest convergence rate of the predictive distribution, here reflected in quite narrow credible intervals. In fact, a (too) fast predictive convergence may underlie a learning mechanism that does not give enough weight to the information in the data:
  $G^*$ is not included in the credible intervals in panels (a) and (b). Frequentist coverage is known to be a delicate issue in Bayesian nonparametric inference (see for example \cite{szabo2015} and the related discussion), but here we can give novel insights, that we find quite intriguing. 
In fact, for c.i.d. as well as for exchangeable data, the credible intervals express the speed of convergence of the predictive distribution (consider expression (\ref{eq:jointPosterior})).  A fast convergence of $G_n(\cdot)(x_{1:n})$ (explicitly denoting the dependence on the data) means that, 
given $x_{1:n}$, there is little uncertainty on the limit $G(\cdot)(x_1, x_2, \ldots)$; therefore, one has narrow credible intervals. 
Clearly, if the predictive distribution is not very sensitive to the data, it will be more stable and generally converge more rapidly; but this may imply a poor learning mechanism, leading to credible intervals that may fail to properly quantify the uncertainty when the data are i.i.d. according to a true distribution. This behavior is evident for Newton's predictive rule $G_n$, because its recursive form clearly outlines the weight given to the current observation as expressed by $\alpha_n$; but the same predictive properties and the same interpretation of credible intervals hold for exchangeable data. Thus, a proper balance between the learning and the predictive convergence rates is also a crucial issue in Bayesian inference. 
Slowly decaying $\alpha_n$ (panels (c) and (d) in  Figure \ref{fig:Gn1000}) give more weight to the current information in the predictive distribution $G_n$; this implies more uncertainty around its limit, thus a larger asymptotic variance and wider credible intervals.

We observed the same behavior for different mixing distributions $g^*$ and varying values of $n$. Figure \ref{fig:GnMultimodale} shows the recursive estimate $G_n$ and the asymptotic $95\%$ marginal credible intervals  for the same simulation setting as in Figure \ref{fig:multimodale-s1}. Here, $g^*$ is multimodal and the sample size is $n=5000$. As in Figure \ref{fig:Gn1000}, the credible intervals in panels (a) and (b) are narrow and do not include the true mixing distribution (dashed curve). 

%%%% G bimodale, n=10mila
%\begin{figure*}[t!]  
% \centering 
%\begin{subfigure}[t]{0.45\textwidth} \label{fig:order-bimodal-10mila}
%	\includegraphics[width=\linewidth]{Rplot14ott-dens10mila-30perm.pdf}
%	\caption{Solid curve: Mixing density estimate (mean over  30 random permutations, 	plotted in gray). Dashed: true mixing density.}
%\end{subfigure}
%~
%\begin{subfigure}[t]{0.45\textwidth} \label{fig:Gn10mila}
%  \includegraphics[width=\linewidth]{Rplot14ott-G10mila-30perm.pdf}
%  \caption{\centering 
%Estimated mixing distribution $G_n$ (solid curve) with $95\%$  marginal  asymptotic credible intervals for $G(t_j)$, for $t_j$ as evidenced in the plot. Dashed curve: true mixing distribution.}
%\end{subfigure}
% \caption{Location mixture of Gaussians; $\sigma^2=1; n=10,000$. 
% Weights $\alpha_n$ as in Figure \ref{fig:bimodal}.} %\ref{fig:order-bimodal}.}
% \label{fig:Norm10mila}
%\end{figure*}

\section{Further statistical applications and extensions} \label{sec:moreExamples}

We have shown how the recursive algorithm can be framed in a rigorous statistical setting, by reading it as a probabilistic predictive rule. This paves the way to further statistical applications and extensions in various directions. 
In this section, we consider the case where (\ref{eq:newtonAsPred}) fully specifies the law of the process $((X_n, \theta_n))$, by further assuming that $P(\theta_{n+1} \in \cdot \mid x_{1:n}, \theta_{1:n})=P(\theta_{n+1} \in \cdot \mid x_{1:n})$. 

\subsection{Unknown common parameters}

The original version of Newton's algorithm does not envisage unknown common parameters in the mixture's kernels. Extensions  for some specific cases are found in \cite{martinGhosh2008}, and a more systematic proposal is given by \cite{martinTokdar2011}.  However, they do not have the probabilistic model underlying Newton's algorithm, thus the  proposed methods are somehow heuristic, not being  based on a genuine likelihood. 
On the contrary, we can easily extend our probabilistic model (\ref{eq:newtonAsPred}) and obtain proper inference. Let 
$$X_i \mid \theta_i, \xi \indsim f(x \mid \theta_i, \xi)$$
where $\xi$ is an common unknown  parameter. Then Newton's rule assigns the  conditional law $ P(\theta_{n+1} \in \cdot \mid x_{1:n}, \xi)= G_n(\cdot \mid \xi)$, 
where $G_n$ is computed according to the rule (\ref{eq:newton}), with the notation here underlining the dependence on $\xi$. These assumptions imply that the conditional
%a joint 
density of $((\theta_1, X_1), \ldots, (\theta_n, X_n))$,
% conditionally on 
given $\xi$, 
%given by  
is
$ 
p(\theta_1, x_1, \ldots, \theta_n, x_n \mid \xi)= 
\prod_{k=1}^n g_{k-1}(\theta_k \mid \xi) f(x_k \mid \theta_k, \xi), 
$
from which one can obtain the marginal likelihood 
$$ m(x_{1:n} \mid \xi) = \prod_{k=1}^n m_k(x_k \mid \xi, x_{1:k-1}), $$
where $m_k(x_k \mid \xi, x_{1:k-1}) \equiv \int f(x_k \mid \theta_k, \xi) g_{k-1}(\theta_k \mid \xi) d\lambda(\theta_k)$. 
Now, one can naturally derive an empirical Bayes estimator of $\xi$ by maximum marginal likelihood, or proceed in a Bayesian approach by assigning a prior distribution to $\xi$ and computing the corresponding posterior law. This gives a probabilistic basis for  the methods proposed by \cite{martinTokdar2011}.

\subsection{Multiple shrinkage estimation with streaming data}

Estimating the individual parameters $\theta_i$ is another problem of interest. 
Our statistical formulation of the recursive rule allows to obtain the posterior distribution of $\theta_{1:n}$. For the time-dependent mixture model (\ref{eq:modLatentGn})-(\ref{eq:tG2}), this provides exact Bayesian inference on $\theta_i$ (filtering), recursively updated as new data become available. 
In the static setting, one obtains quasi-Bayes inference for the $\theta_i$. Again, the advantage is to allow fast recursive computations with streaming data.

%In mixture models, interest may also be on the estimation of the individual parameters $\theta_i$. When the c.i.d. model is used as a quasi-Bayes approximation of an exchangeable mixture model,   the lack of (approximate) exchangeability in the initial stage of the process must be considered. 
%sensitivity to permutations, \red{and the close approximation to a DP prior when the $X_i$ are very informative about the corresponding $\theta_i$.} 
%Again, the advantage 

The joint posterior density of $\theta_{1:n} $, given $\xi$ and $x_{1:n}$,  is easily obtained as 
\begin{equation} \label{eq:posteriorTheta}
p(\theta_{1:n} \mid \xi, x_{1:n}) = 
\prod_{k=1}^{n}  \frac{f(x_k \mid \theta_k, \xi) g_{k-1}(\theta_k \mid \xi)}{m_k(x_k  \mid \xi, x_{1:k-1})}. 
\end{equation}
In the posterior distribution, one recursively estimates $g$ and uses the estimate as the prior for the new $\theta_k$, independently over the $\theta$'s. 
This is a sort of temporal empirical-Bayes procedure: at time $k$, the sample $x_{1:k-1}$ is used to estimate the \lq \lq prior distribution'' of $\theta_k$ (the latent distribution $\tG_k$ in the time-varying model (\ref{eq:modLatentGn})-(\ref{eq:tG2})); the estimate  $G_{k-1}$ is then used, in an empirical-Bayes fashion, as the prior law for inference on $\theta_k$ based on $x_k$. 
When the common parameter $\xi$ is unknown, inference on the $\theta_i$ can be  solved by plugging the marginal maximum likelihood estimator $\hat{\xi}_n$ into    (\ref{eq:posteriorTheta}) or, in a Bayesian approach, by assigning a prior law on $\xi$ and integrating (\ref{eq:posteriorTheta}) with respect to the posterior distribution of $\xi$. 

\subsection{Multivariate parameters} 

A known limitation of Newton's algorithm is that it requires to evaluate an integral at each step. This can be solved by numerical methods but becomes demanding in the case of multivariate  $\theta$. An interesting class of predictive recursive algorithms that avoid the integral computations has been recently proposed by \cite{hahnMartinWalker2018}, %\red{for density estimation -- VEDERE},
 and an application in a multivariate setting is found in \cite{cappelloWalker2018}. 
Our probabilistic setting can be also exploited for suggesting new computational strategies with multivariate parameters. Here we sketch a simple  Monte Carlo scheme. Although we do not expand further nor evaluate the Monte Carlo error, simulation results %(not shown  \red{ho messo una figura \ref{fig:Bin-MCvsNewton})} 
are encouraging, showing very good approximations.

%\begin{figure} \label{fig:Bin-MCvsNewton}
%\includegraphics[scale=1]{RplotBin-MCvsNewton}
%\caption{Binomial data: $X_i \mid \theta_i \indsim$ Binomial($M, \theta_i$), with $M=20$, sample size $n=1000$. Recursive estimate $g_n$ of the mixing density (Monte Carlo (solid) and numerical integration (dashed)). Split-sample weights $\alpha_n$, with $M=600, \alpha=100$.}
%\end{figure}
%

\begin{figure}[h] 
%% BINOMIALE -- Dati salvati come xBin e ThetaBin
\centering
\includegraphics[width=.7\textwidth,height=6cm]{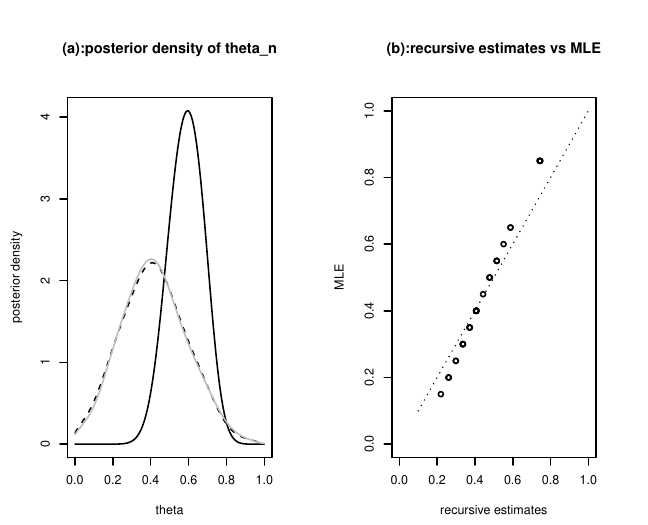}
 \caption{Binomial data: $X_i \mid \theta_i \indsim$ Binomial($M, \theta_i$), $M=20$; sample size $n=1000$. Split-sample weights $\alpha_n$ with $N=500, \gamma=3/4$; $\alpha=100$. 
Panel (a) Posterior density of $\theta_n \mid x_{1:n}$ (solid) and predictive density of $\theta_n \mid x_{1:n-1}$ (Monte Carlo (dashed) and numerical (gray) integration). The two predictive densities are almost overlapping.
 Panel (b): Recursive estimates of the last $30$ parameters $\theta_i$ versus the maximum likelihood estimates.}
\label{fig:binomial} 
 \end{figure}

Notice that one can write the recursive rule (\ref{eq:recursiveDensities}) as
$$ 
g_n(\theta ) = g_{n-1}(\theta) \, \left( 1+ \alpha_n (\frac{f(x_n \mid \theta)}{m_n(x_n \mid x_{1:n-1})} - 1)\right).
$$ 
%where for brevity we here drop the possible dependence on common parameters $\xi$. 
Iterating, one gets
\begin{equation} \label{eq:newtonMC}
g_n(\theta ) = g_0(\theta) \, \prod_{k=1}^n \left( 1+ \alpha_k (\frac{f(x_k \mid \theta)}{m_k(x_k \mid x_{1:k-1})} - 1)\right) , 
\end{equation}
and can obtain the following recursions for the integrals $m_k(x_k \mid x_{1:k-1})$:
\begin{eqnarray*}
m_1(x_1) &=& \int f(x_1 \mid \theta) g_0(\theta) d \lambda(\theta), \\
m_k(x_k \mid x_{1:k-1}) &=& \int  f(x_k \mid \theta)\prod_{i=1}^{k-1} \left( 1+ \alpha_i (\frac{f(x_i \mid \theta)}{m_i(x_i \mid x_{1:i-1})} - 1)\right) g_0(\theta) d\lambda(\theta), \quad k \geq 2.  
\end{eqnarray*}

 One can then envisage a Monte Carlo scheme for recursively approximating the integrals $m_k$ and for sampling from the posterior distribution. 
 It is indeed enough to sample from the prior density $g_0$. Let $(\theta_1^*, \ldots, \theta_M^*)$ be a pseudo-random sample from $g_0$. Then a Monte Carlo estimate of the integrals $m_k(x_k \mid x_{1:k-1})$ can be obtained starting from 
$ \hat{m}_1(x_1)= \sum_{j=1}^M f(x_1 \mid \theta_j^*)/M$ and 
and {\em recursively} computing
$$ \hat{m}_k(x_k\mid x_{1:k-1}) = 
\frac{\sum_{j=1}^M f(x_k \mid \theta_j^*)\prod_{i=1}^{k-1} \left( 1+ \alpha_i (\frac{f(x_i \mid \theta_j^*)}{\hat{m}_i(x_i \mid x_{1:k-1})} - 1)\right)}{M}
$$
for $k > 1$. 
The Monte Carlo estimates $\hat{m}_k(x_k \mid x_{1:k-1})$ are fairly easily computed even for multivariate parameters $\theta_i$, and can be used for recursively evaluating $g_k(\theta)$. A Monte Carlo evaluation of the posterior density of $\theta_k$ can also be recursively obtained as  
$$ \hat{p}_{g_{k-1}}(\theta \mid x_k) = \frac{g_{k-1}(\theta) f(x_k \mid \theta)}{\hat{m}_k(x_k \mid x_{1:k-1})},
$$
and using expression (\ref{eq:newtonMC}) one can envisage sampling from the posterior distribution of $\theta_k$ by sampling from $g_0$.

We illustrate the procedure for a small example. The data are generated as 
$X_i \mid \theta_i \indsim$ Binomial($M, \theta_i$), $i=1, \ldots, n$, with $M=20$ and $n=1000$; the $\theta_i$ are i.i.d. from a Beta distribution with parameters $(3,4)$. 
The initial distribution $G_0$ is Uniform($0,1$). 
Interest is in recursive estimation of the last $K$ values of $\theta$. We use the recursive rule $g_n$ with weights $\alpha_n$ fixed by split-sample with $N=500, \gamma=3/4$, and $\alpha=100$. The results are shown in Figure \ref{fig:binomial}. 
Panel (a) shows the posterior density of $\theta_n$, given $x_{1:n}$, together with the predictive density $g_n$ of $\theta_n$ given $x_{1:n-1}$, evaluated through Monte Carlo integration (dashed; Monte Carlo sample of size $100,000$) and through numerical integration (gray).  The two predictive densities are almost overlapping. 
In panel (b), we plot
%Subfigure (\ref{fig:binomial}b) plots 
the recursive estimates $E(\theta_i \mid x_{1:i})$ against the maximum likelihood estimates (MLE) $\hat{\theta}_i=x_i/M$, for $i=n-30+1, n$. The shrinkage effect is evident. Working  with simulated data, we have the true values of the $\theta_i$ and can compute the MSE, which is  $0.01062$ for the MLE and 
$0.00648 $ for the recursive estimates. 
% MSElik; MSEnewton
% [1] 0.01061759
% [1] 0.006459924}

\section{Discussion}
\label{sec:discussion}

Due to its simplicity and good practical performance, Newton's algorithm is quite popularly used in problems involving hidden variables. We have proposed a novel approach that develops the algorithm into a quasi-Bayes method, and makes the user aware of the modeling assumptions implicitly made. We believe that our approach can also be useful in other settings. 

Explicit results on the probability law of the asymptotic mixing distribution $G$, although  difficult to obtain, would give a more complete description of the prior implied by the recursive predictive rule, and the construction could be further extended to characterize 
 novel priors on the space of  absolutely continuous  distributions, for Bayesian nonparametrics. Modifications of the algorithm could be envisaged, for example by initializing the procedure with exact computations from the DP mixture model, in order to control the prior distribution on $G$.  

The lack of exchangeability of Newton's algorithm has been addressed in the literature by taking an average of the recursive estimates over a number of random permutations of the original data, although this procedure sacrifices the recursive nature of the computations. Our approach may be useful to interpret this modified algorithm in a proper statistical framework. We can formalize the permutation-based algorithm as defining a new predictive rule, that assumes $\theta_1 \sim G_0$ and, for any $n \geq 1$,
$
\theta_{n+1} \mid x_{1:n} \sim \bar{G}_n(\cdot)= {\sum_{\pi} G_n^{\pi}(\cdot)}/{M} , 
$
where $G_n^{\pi}$ is the estimate (\ref{eq:newton})  obtained for a random permutation $\pi$ of $x_{1:n}$ and $M$ is the total number of permutations considered. 
Adding for simplicity the assumption that $\theta_{n+1}$ is conditionally independent on $\theta_{1:n}$ given $x_{1:n}$, by Ionescu-Tulcea Theorem this predictive rule characterizes a new probability law for the process $((X_n, \theta_n))_{n \geq 1}$.  One may thus follow our predictive approach to study this  new process and 
develop the permutation-based modification of Newton's algorithm into a proper statistical method, 
that may interestingly define a new prior.
% We also conjecture that one could use our techniques to obtain an  asymptotic approximation of the posterior distribution of $G(t)$ given $x_{1:n}$ for this new model, as a Gaussian centered at  $\bar{G_n}(t)$.
%%This would develop the permutation-based modification of Newton's algorithm into a proper statistical method, with uncertainty quantification through credible intervals. Our results would still be useful to give an advise on the choice of $M$. 

%%Another direction of research suggested by our work refers to frequentist coverage properties.  
%The remarks in Section \ref{sec:inferenceG} on the behavior of the quasi-Bayes asymptotic credible intervals also holds in exchangeable setting. This appears to suggests a novel perspective on frequentist coverage properties of Bayesian procedures, where the challenge is to have the right balance between the learning rate and the predictive convergence rate.

A computational limitation of Newton's algorithm is that it requires to evaluate an integral at each step. We have described a simple Monte Carlo approximation, and plan to further explore this issue in future work. 
Extensions of  our study to the class of algorithms proposed by \cite{hahnMartinWalker2018}, as well as developments for multivariate mixtures and dependent mixture models, possibly exploiting theoretical results on {\em partially c.i.d.} sequences (\cite{fortiniPetronePolina2017}) present interesting direction for future research.

\bigskip

\noindent {\bf Acknowledgments}  We sincerely thank the AE and the reviewers for their stimulating
%\red{cut? thoughtful and constructive} 
comments. 
%that lead to an improved version of our manuscript?. 
We also express our gratitude to Lorenzo Cappello and Stephen Walker for many inspiring discussions and thank Isadora Antoniano Villalobos for her helpful comments.
%This work benefited from discussions with Lorenzo Cappello and Stephen Walker, to whom we express our thanks. 
The authors acknowledge financial support by the PRIN grant 2015SNS29B.

\bibliographystyle{plainnat}   %may change style
\bibliography{pcid}
 
\newpage
\pagenumbering{arabic}
\setcounter{section}{0}
\setcounter{equation}{0}
\renewcommand{\thesection}{A\arabic{section}}
\renewcommand{\theequation}{A\arabic{equation}}
\vspace{3cm}

	\large{
\bf \noindent Appendix - Quasi Bayes Properties of a procedure for sequential learning \\in mixture models}\\
\\
Sandra Fortini and Sonia Petrone\\
{\sl Bocconi University, Milan, Italy}\\
sonia.petrone@unibocconi.it

\vspace{1cm}

\section{Technical details and proofs} \label{sec:proofs}
%This Appendix contains proofs and technical details. 
We will make use of the following notion (\cite{BertiPratelliRigo2004LimThForIID}).

\begin{defn}
A sequence of random variables $(X_n)$ is {\em conditionally identically distributed with respect to the filtration $\cF$} ($\cF$-c.i.d.) if it is adapted to $\cF$ and
$$ E[ h(X_{n+k}) \mid \mathcal F_n] = E[ h(X_{n+1}) \mid \mathcal F_n], 
$$
for all $k \geq 1$, $n \geq 0$ and all bounded measurable functions $h: \mathbb X \rightarrow \mathbb{R}$. 
\end{defn}

When $\cF $ is the natural filtration of $(X_n)$, the sequence is said to be c.i.d. An $\cF$-c.i.d. sequence is also c.i.d. 
Unless otherwise specified, in the sequel we denote by $\cF$ the natural filtration of $(X_n)$.
%denote by $\cF_0$ the trivial sigma field, by $\cF_n$ the sigma field generated by $X_{1:n}$ ($n\geq 1$), and let  $\cF_\infty=\vee_n \cF_n$. 

\medskip

%\begin{defn} 
%\label{def:cid}
%Let $\cF= (\cF_n)_{n \geq 1}$ be a filtration. A sequence of random variables $(X_n)$ is {\em conditionally identically distributed} with respect to the filtration $\cF$ ($\cF$-c.i.d.) if it is adapted to $\cF$ and
%$$ E[ h(X_{n+k}) \mid \mathcal F_n] = E[ h(X_{n+1}) \mid \mathcal F_n], 
%$$
%for all $k \geq 1$, $n \geq 0$ and all bounded measurable functions $h: \mathbb X \rightarrow \mathbb{R}$ 
%\end{defn} 
%Less formally, a sequence $(X_n)$ adapted to $\cF$ is $\cF$-c.i.d. if the $X_i$ are identically distributed and $X_{n+k} \mid \cF_n \eqdist X_{n+1} \mid \cF_n, \quad k \geq 1, n \geq 1$.  
%When $\cF $ is the natural filtration of $(X_n)$, the sequence is said to be c.i.d. An $\cF$-c.i.d. sequence is also c.i.d.  

\noindent {\sc Proof of Theorem \ref{th:urnmart}}. 

(i) For every $A\in\mathcal B(\Theta)$ and every $n\geq 0$,
%\red{one has $E(G_1(A)=G_0(A)$} and for every \red{$n\geq 1$}
\begin{equation} \label{eq:Gn-martingale}   
E(G_{n+1}(A)\mid \cF_n)=
(1-\alpha_{n+1})G_n(A)+\alpha_{n+1} E(P(\theta_{n+1}\in A\mid \cF_{n+1})\mid  \cF_n)=G_n(A).
\end{equation}
Hence, the sequence $(G_n)$ is a measure valued martingale, under $P$, with respect to the natural filtration of $(X_n)$. %$(\cF_n)$.
By Lemma 7.14 in \cite{aldous1985}, there exists a random probability measure $G$ such that $G_n$ converges $P$-a.s. to $G$, in the topology of weak convergence. 

 (ii) Since, for every $A$, $(G_n(A))$ is uniformly bounded, it is a closed martingale. Thus, for every $n\geq 0$ and $k\geq 1$,
\[
P(\theta_{n+k}\in A\mid \cF_n)=E(P(\theta_{n+k}\in A\mid
\cF_{n+k-1}) \mid \cF_{n})=E(G_{n+k-1}(A)\mid \cF_{n})=
E(G(A)\mid \cF_{n}).
\]
$\square$

\bigskip

\noindent {\sc Proof of Proposition \ref{prop:hconv}}. 

Let  $Z$ be a random variable such that $Z \mid \cF_\infty \sim G$. We have
\begin{eqnarray*}
E\left( \int_\Theta h(z) dG(z) \left| \cF_n\right.\right)&=& E\left( E(h(Z) \mid \cF_\infty) \mid \cF_n \right) = 
E(h(Z) \mid \cF_n) \\
&=& \int_\Theta h(z) P(Z\in dz \mid \cF_n), 
\end{eqnarray*}
where the last equality follows from the fundamental property of regular conditional distributions (see e.g. \cite{aldous1985}, eq.(2.4)). Noticing that  
$P(Z \in \cdot \mid \cF_n)= E( P(Z  \in \cdot \mid \cF_\infty) \mid \cF_n)= E(G(\cdot) \mid \cF_n)$, 
we obtain 
\begin{equation} \label{eq:expectationGn}
\int_\Theta h(\theta) dG_n(\theta) = E\left(\int_\Theta h(\theta)dG(\theta) \mid  \mathcal F_{n}\right). 
\end{equation}
Since 
$$
E\left(\int_\Theta h(\theta)dG(\theta) \mid  \mathcal F_{n}\right) \rightarrow 
E\left(\int_\Theta h(\theta)dG(\theta) \mid  \mathcal F_\infty \right)\quad \mbox{$P$-a.s.},
$$
and since $G$ is $\cF$-measurable, then $\int_\Theta h(\theta) dG_n(\theta)\rightarrow  \int_{\Theta} h(\theta) d G(\theta)$ $P$-a.s.

To prove the last assertion, notice that, from (\ref{eq:expectationGn}), 
$\int_\Theta |h(\theta)| dG_0(\theta)= E(\int_\Theta |h(\theta)|  dG(\theta))$. 
Thus, if $\int_\Theta |h(\theta)| dG_0(\theta) < \infty$, then the non-negative quantity $\int_\Theta |h(\theta)| dG(\theta)$ is $P$-a.s. finite. \\
$\square$

\bigskip

\noindent {\sc Proof of Theorem \ref{th:cidconv}}

(i) To prove that $(X_n)$ is c.i.d., it is enough to show that 
$P(X_{n+2} \in B \mid \cF_{n})=P(X_{n+1} \in B \mid \cF_{n})$, for any $n \geq 0$ and any $B$. This is a consequence of the conditional independence of the $X_i$, given $(\theta_n)$, and of Theorem \ref{th:urnmart}. Indeed, denoting by $F(\cdot \mid \theta)$ the distribution corresponding to the density $f(\cdot \mid \theta)$, 
we have 
	\begin{eqnarray*}
P(X_{n+2} \in B \mid \cF_{n}) &=&  E(P(X_{n+2} \in B \mid \theta_{n+2}, \cF_{n}) \mid \cF_{n})
=\int_\Theta F(B \mid \theta) P(\theta_{n+2}\in d\theta \mid \cF_{n}) \\ 
&=& \int_\Theta F(B \mid \theta) P(\theta_{n+1}\in d\theta \mid \cF_{n}) = P(X_{n+1} \in B \mid \cF_{n}), 
\end{eqnarray*}
where the third equality follows from (\ref{eq:Gn-martingale}). 

(ii) Let $(t_j, j \in J)$  be a countable dense set of points in $\mathbb X$. 
By Proposition \ref{prop:hconv}, $ F_{G_n}(t_j) \rightarrow F_G(t_j)$,  
for every $\omega \in \Omega_j$ with $P(\Omega_j)=1$. Now, let $\Omega^*=\cap_j \Omega_j$. Being $J$ countable, $P(\Omega^*)=1$,  and for any $\omega \in \Omega^*$, $F_{G_n}(t_j) \rightarrow F_G(t_j)$ for all $t_j$. 
For distribution functions, convergence on a countable dense set implies weak convergence. Therefore, we have that, $P$-a.s., $F_{G_n}$ converges to $F_G$ in the topology of weak convergence.  Now, notice that $F_G$ is $P$-a.s. absolutely continuous, with density $f_G$. By Theorem 1 in \cite{BertiPratelliRigo2013}, $P$-a.s. weak convergence of the predictive measures  to an absolutely continuous random measure  implies that the convergence also holds in total variation. Therefore, $P$-a.s., $F_{G_n}$ converges to $F_G$ in total variation, which is equivalent to $f_{G_n} {\overset{L_1} \rightarrow} f_G$. 

\noindent (iii) Convergence of the predictive distributions to the random probability measure  $F_G$ implies that $(X_n)$ is asymptotically exchangeable, with directing measure $F_G$ (\cite{aldous1985}, Lemma 8.2). \\
$\square$

\bigskip

\noindent {\sc Proof of Theorem \ref{th:asympttheta}}

Let $\mathcal H$ be a a countable, convergence determining class   of bounded continuous functions, and $k$ a positive integer. 
By Theorem \ref{th:cidconv},  $(X_n)$ is asymptotically exchangeable with directing random measure $F_G$; therefore, for $P$-almost all $\omega=(x_1,x_2,\dots)$
\begin{equation}
\label{eq:exchangX}
E(\prod_{i=1}^kh_i(X_{n+i})\mid x_{1:n})\rightarrow \prod_{i=1}^k
 \left(\int_\Theta  \int_{\mathbb X}h_i(z_i)f(z_i\mid s_i)d\mu(z_i) dG(s_i)(\omega)\right), 
\end{equation}
for every $h_i\in\mathcal H$.
Let $\omega=(x_1,x_2,\dots)$ be fixed in such a way that the above equations hold.   
Then, for every $j=1,\dots,k$, the sequence of probability measures $\left(P( \theta_{n+j}\in\cdot \mid x_{1:n})\right)$ is tight. It follows that the sequence of joint conditional distributions $\left(P( \theta_{n+1:n+k}\in\cdot \mid x_{1:n})\right)$ is tight.
Thus, for every increasing sequence of integers, there exists a subsequence $(n_j)$ and a probability measure $Q(\omega)$ such that  
$ %\begin{equation} \label{eq:conv1asyex}
P(( \theta_{n_j+1},\dots, \theta_{n_j+k})\in \cdot\mid  x_{1:n_j})\rightarrow   Q(\omega)$.  
The proof is complete if we can show that 
\begin{equation} \label{eq:conv4asyex}
Q(A_1\times \dots \times  A_k)(\omega)=\prod_{i=1}^k G(A_i)(\omega), \quad \mbox{for every $A_1,\dots A_k$}, 
\end{equation}
because (\ref{eq:conv4asyex}) implies that the conditional law of $( \theta_{n_j+1},\dots, \theta_{n_j+k})$ converges weakly to the product measure $G^k$. The sequence of random variables 
$
\left(\prod_{i=1}^k \int_{\mathbb X}h_i(z_i)f(z_i\mid\theta_{n+i})d\mu(z_i)\right)$
is uniformly bounded and, therefore, it is uniformly integrable. Hence, 
$$ % \begin{equation} \label{eq:conv2asyex}
\begin{aligned}
\int_{\Theta^k}\prod_{i=1}^k\int_{\mathbb X} h_i(z_i) f(z_i\mid s_i) d\mu(z_i) &  P(( \theta_{n_j+1},\dots \theta_{n_j+k})\in(ds_1,\dots,ds_k)\mid x_{1:n_j})
\\ &\rightarrow \int_{\Theta^k} \prod_{i=1}^k \int_{\mathbb X} h_i(z_i) f(z_i\mid s_i) d\mu(z_i)  dQ(s_1,\dots,s_k)(\omega),
\end{aligned}
$$
for every $h_1,\dots,h_k\in \mathcal H$.
On the other hand, by (\ref{eq:exchangX}),
$$
\begin{aligned}
\int_{\Theta^k}\prod_{i=1}^k\int_{\mathbb X}h_i(z_i)f(z_i\mid s_i)d\mu(z_i)&P(( \theta_{n_j+1},\dots \theta_{n_j+k})\in(ds_1,\dots,ds_k)\mid x_{1:n_j})
\\&\rightarrow \int_{\Theta^k}\prod_{i=1}^k\int_{\mathbb X}h_i(z_i)f(z_i\mid s_i)d\mu(z_i) dG(s_1)(\omega) \dots dG(s_k)(\omega).
\end{aligned}
$$
Hence, for every $h_1,\dots,h_k\in\mathcal H$,
\begin{equation} \label{eq:conv3asyex}
\begin{aligned}
\int_{\Theta^k}\prod_{i=1}^k\int_{\mathbb X}h_i(z_i)f(z_i\mid s_i)  & d\mu(z_i) dQ(s_1,\dots,s_k)(\omega)\\
&=\int_{\Theta^k}\prod_{i=1}^k\int_{\mathbb X} h_i(z_i)f(z_i\mid s_i)d\mu(z_i)dG(s_1)(\omega) \dots dG(s_k)(\omega).
\end{aligned}
\end{equation}
Since the model is identifiable,
the class 
\[
\left\{\int h(z)f(z\mid  \theta )d\mu(z) : h\in \mathcal H    \right\}
\]
is separating for $\mathcal P(\Theta)$. It follows that the class 
\[
\left\{\prod_{i=1}^k\int h_i(z)f(z\mid \theta)d\mu(z):h_i\in \mathcal H,i=1,\dots,k    \right\}
\]
is separating for $\mathcal P(\Theta^k)$ (\cite{ethierKurtz1986}, Proposition 3.4.6).
Thus, (\ref{eq:conv3asyex}) implies (\ref{eq:conv4asyex}).\\
$\square$

\bigskip

The proof of Theorem \ref{th:absolcont} is based on the following Lemmas \ref{th:appendix1} and \ref{th:appendix2},
%the following Lemmas,
 which are extensions  of Theorems 1 and 4 in  \cite{BertiPratelliRigo2013}.
 %, the main difference being that we do not assume that the sequence of random variables involved is c.i.d.  
%\red{(a property that does not hold, in general, for the sequence $(\theta_n)$ in Newton's model).} .
%If $G_0$ is absolutely continuous with respect to a sigma-finite measure $\lambda$ on $\Theta$,  denoted $G_0 \ll \lambda$, then $G_n \ll \lambda$, and the corresponding density $g_n$ satisfies Newton's recursive rule (\ref{eq:recursiveDensities}). 
%It is easy to verify that, for any fixed $\theta$,  the sequence $(g_n(\theta))$ is a martingale \red{under the c.i.d. law $P$}.  Since $g_n(\theta)$ is non-negative, there exists a function $g^*(\theta)$ such that, for every $\theta$,  $g_n(\theta)$ converges to $g^*(\theta)$ a.s. However, this fact is not sufficient to conclude  that $G \ll \lambda$. 
%Extending a remarkable result by \cite{BertiPratelliRigo2013}, we 

Lemma \ref{th:appendix1} shows that $G \ll \lambda$ requires that $g_n$ converges in $L^1$ or, equivalently, that $G_n$ converges to $G$ in total variation.  
%This is guaranteed in . 
Then, Lemma \ref{th:appendix2} gives sufficient conditions for $G_n \rightarrow G$ in total variation. % and $G\ll\lambda$, a.s.
%The two lemmas are extensions of Theorems 1 and 4 in  \cite{BertiPratelliRigo2013}, the main difference being that we do not assume that the sequence of random variables involved is c.i.d.  
%\red{(a property that does not hold, in general, for the sequence $(\theta_n)$ in Newton's model).} 
%\red{CUT? To be c.i.d. with respect to a filtration  $\cF$, a sequence has to be adapted to $\cF$, which, in general, is not true for $(\theta_n)$ in Newton's model and the filtration of interest.}
The proofs of the lemmas can be fairly easily obtained, adapting the ones of Theorems 1 and 4 in  \cite{BertiPratelliRigo2013} and directly requiring the martingale property of the sequence of random measures $Q_n$, which is otherwise implied by the $\cF$-c.i.d. property. 
%However, the proofs of Theorems 1 and 4 in \cite{BertiPratelliRigo2013} can be adapted,
% by directly requiring the martingale property of the sequence of the random measures $Q_n$, which is otherwise implied by the $\cF$-c.i.d. property.
Thus, the following two lemmas are provided without additional proof.

\begin{lemma}
	\label{th:appendix1}
	Let $\lambda$ be a sigma-finite measure on a Polish space $S$. For any $n$, let $Q_n$ be a random measure on $S$ such that the sequence $(Q_n)$ is a measure-valued martingale, under $P$, with respect to a filtration $(\mathcal F_n)$, and let $Q$ be its limit. 
	Then $Q \ll \lambda$, $P$-a.s.  if and only if, $P$-a.s.,   
	$Q_n \ll \lambda$ for every $n$ and $Q_n$ converges to $Q$ in total variation.
\end{lemma}

\begin{lemma}
	\label{th:appendix2}
	Let $\lambda, Q_n$ and $Q$ be as in Lemma \ref{th:appendix1}.
	Assume that $Q_n \ll \lambda$, $P$-a.s., for every $n$, with density  $q_n$. 
	Then $Q\ll \lambda$, $P$-a.s.,  if and only if for every compact $K$ such that $\lambda(K)<\infty$, 
	$q_n$ is, $P$-a.s., a  function on $S$ uniformly integrable with respect to $\lambda_K$, 
	where $\lambda_K(\cdot)=\lambda(\cdot\cap K)$ is the restriction of $\lambda$ on $K$. \\
	In particular, $Q\ll\lambda$, $P$-a.s., if, for every $K$ compact, there exists $p>1$ such that, $P$-a.s., 
	\begin{equation} \label{eq:appendix2}
		\sup_n\int_K q_n(x)^pd\lambda(x)<\infty. 
	\end{equation}
	A sufficient condition for (\ref{eq:appendix2}) is 
	$$
	\sup_nE(\int_K q_n(x)^pd\lambda(x))<\infty.
	$$
\end{lemma}

\noindent{\sc Proof of Theorem \ref{th:absolcont}.} 

The thesis follows from Lemmas \ref{th:appendix1} and \ref{th:appendix2}, if we can show that
\begin{equation}
\label{eq:suffcond}
\sup_nE\left(\int_Kg_n(\theta)^2d\lambda(\theta)\right)<\infty\quad\mbox{ for every }K\mbox{ compact, satisfying }\lambda(K)<\infty.
\end{equation}
Let $K$ be a fixed compact set, with $\lambda(K)<\infty$. It holds
\[
E\left( \int_Kg_n(\theta)^2d\lambda(\theta)  \right)   =\int_K E(g_n(\theta)^2)d\lambda(\theta)=\int_KE( E(g_n(\theta)^2\mid \cF_{n-1}))d\lambda(\theta).
\]
By the martingale property of the sequence $(g_n)$, and Jensen inequality, we obtain 
\[
\begin{aligned}
E\left(g_n(\theta)^2\mid X_{1:n-1} \right)
&= g_{n-1}(\theta)^2 
E\left( \left[1+\alpha_n\left(\frac{f(X_n\mid\theta)}{\int_\Theta f(X_n\mid\theta') g_{n-1}(\theta') d\lambda(\theta')}-1\right)\right]^2  \mid \cF_{n-1} \right)\\
%&=g_{n-1}(\theta)^2\left[ 1+\alpha_n^2E\left(\left(\frac{f(X_n\mid\theta)}{\int_\Theta f(X_n\mid\theta')g_{n-1}(\theta')d\lambda(\theta')}-1\right)^2  \mid X_{1:n-1} \right) \right]\\
%&\leq g_{n-1}(\theta)^2\left[ 1+\alpha_n^2\left(1+\int_{\mathbb X}  \frac{f(x\mid\theta)^2}{\int_\Theta f(x\mid\theta')g_{n-1}(\theta')d\lambda(\theta')} d\mu(x)  \right)\right]\\
&\leq g_{n-1}(\theta)^2\left[ 1+\alpha_n^2\left(1+\int_{\mathbb X}  \frac{f(x\mid\theta)^2}{\int_K f(x\mid\theta')g_{n-1}(\theta')d\lambda(\theta')} d\mu(x)  \right)\right]\\
&\leq g_{n-1}(\theta)^2\left[ 1+\alpha_n^2\left(1+\int_{\mathbb X} \int_K \frac{f(x\mid\theta)^2}{ f(x\mid\theta')}g_{n-1}(\theta')d\lambda(\theta') d\mu(x)  \right)\right]\\
&\leq g_{n-1}(\theta)^2\left[ 1+\alpha_n^2\left(1+\sup_{\theta_1,\theta_2\in K}\int_{\mathbb X} \frac{f(x\mid\theta_1)^2}{ f(x\mid\theta_2)} d\mu(x)  \right)\right].
\end{aligned}
\]
Therefore
\[
E\left(\int_Kg_n(\theta)^2d\lambda(\theta)\right)\leq E\left(\int_Kg_{n-1}(\theta)^2d\lambda(\theta)\right)\left[ 1+\alpha_n^2\left(1+\sup_{\theta_1,\theta_2\in K}\int_{\mathbb X} \frac{f(x\mid\theta_1)^2}{ f(x\mid\theta_2)} d\mu(x)  \right)\right]
\]
Iterating, we obtain
\[
E\left(\int_Kg_n(\theta)^2d\lambda(\theta)\right)\leq \int_Kg_0(\theta)^2d\lambda(\theta)\prod_{i=1}^n\left( 1+\alpha_i^2 M_K\right),
\]	
with $M_K=\left(1+\sup_{\theta_1,\theta_2\in K}\int_{\mathbb X} \frac{f(x\mid\theta_1)^2}{ f(x\mid\theta_2)} d\mu(x)  \right)$, which is finite by the assumption (\ref{eq:suffcont3}). By  (\ref{eq:suffcont1}),  %and (\ref{eq:suffcont2}),
$\sup_nE(\int_Kg_n(\theta)^2d\lambda(\theta))<\infty$.\\
$\square$

\bigskip

\noindent {\sc Proof of Lemma \ref{lemma:V}}. 

Since $G_n(A)\rightarrow G(A)$ $P$-a.s., it remains to show that
$$
\int_{\mathbb X} P_{G_n}(A \mid x)^2 dF_{G_n}(x) \rightarrow 
\int_{\{x:f_G(x)\neq 0\}} P_{G}(A \mid x)^2 dF_{G}(x)\quad \mbox{$P$-a.s.}
$$
By Theorem \ref{th:cidconv}, $F_{G_n}$ converges to $F_G$ in total variation, $P$-a.s. Therefore,
$$
\int_{\{x:f_G(x) = 0\}} P_{G_n}(A \mid x)^2 dF_{G_n}(x) 
\leq F_{G_n}(\{x:f_G(x)=0\}) \rightarrow F_{G}(\{x:f_G(x)=0\})=0  \quad \mbox{$P$-a.s.}
$$
Then, denoting $\mathbb X_0=\{x:f_G(x)\neq 0\}$, we have
\begin{eqnarray*}
	&&\left|\int_{\mathbb X_0} P_{G_n}(A \mid x)^2 dF_{G_n}(x) - \int_{\mathbb X_0} P_{G}(A \mid x)^2  dF_G(x)\right|	\\
	&\leq & 
	\int_{\mathbb X_0} |f_{G_n}(x) - f_G(x)| d\mu(x) + \left|\int_{\mathbb X_0} P_{G_n}(A \mid x)^2(x) dF_{G}(x) - \int_{\mathbb X_0} P_{G}(A \mid x)^2 dF_G(x)\right|. \\
\end{eqnarray*}
The first term converges to zero since $f_{G_n}$ converges to $f_G$ in $L_1$, by Theorem \ref{th:cidconv}. The second term converges to zero by dominated convergence theorem. Thus, the thesis follows.\\
$\square$

\bigskip

\noindent{\sc Proof of Theorem \ref{th:rategeneral}}.

For every $n\geq 1$, let 
$$M_{n,j}=\left\{
\begin{array}{ll}
\sqrt{r_{n}}(G_n(A)-G_{n+j-1}(A))&j\geq 1\\
0&j=0,
\end{array}
\right.
$$
and let 
$$\mathcal F_{n,j}=\left\{
\begin{array}{ll}
\mathcal F_{n+j-1}& j\geq 1\\
\mathcal F_n &j=0.
\end{array}
\right.
$$
For every $n\geq 1$, $(M_{n,j})_{j\geq 0}$ is a zero-mean martingale, under $P$, with respect to the filtration $(\mathcal F_{n,j})_{j\geq 0}$ and $\mathcal F_{n,1}=\mathcal F_n\subset \mathcal F_{n+1}=\mathcal F_{n+1,1}$.
Let
$$
Z_{n,j}\equiv M_{n,j}-M_{n,j-1}\mbox{ for }j\geq 1,\quad U_n\equiv\sum_{j\geq 1}Z_{n,j}^2,\quad Z_n^*\equiv\sup_{j\geq 1}|Z_{n,j}|.
$$
The thesis follows from Theorem A.1 in \cite{crimaldi2009} if we can show that $(Z_n^*)$ is dominated in $L^1$ and that $(U_n)$ converges $P$-a.s. to $V_A$.

%Theorem A.1 in \cite{crimaldi2009} requests to show that $(Z_n^*)$ is dominated in $L^1$ and $(U_n)$ converges a.s. to $V_A$.

By definition,
$Z_{n,1}=0$ and, for $j\geq 2$,
\[
\begin{aligned}
Z_{n,j}&=\sqrt{r_{n}}\left(G_{n+j-2}(A)-G_{n+j-1}(A)\right)\\
&=\sqrt{r_{n}}\alpha_{n+j-1}\left(G_{n+j-2}(A) -\frac{\int_Af(X_{n+j-1}\mid\theta)dG_{n+j-2}(\theta)}{\int_\Theta f(X_{n+j-1}\mid\theta)dG_{n+j-2}(\theta)} \right).
\end{aligned}
\]
Since $\sqrt{r_n}\sup_{k\geq n}\alpha_k\rightarrow 0$, then  $(Z_n^*)$ is dominated in $L^1$.\\
To prove that $(U_n)$ converges $P$-a.s. to $V_A$, we employ Lemma A.1 in \cite{crimaldi2016fluctuation}. To be consistent with the notation therein, let us set $b_1=r_1$ and, for $k\geq 1$,
$$ b_{k+1}=r_{k}\quad\mbox{ and }\quad a_k=\frac{1}{b_{k}^2\alpha_k^2}.
$$
Then, we can write
\[
\begin{aligned}
U_n&=b_{n+1}\sum_{j\geq 2}\alpha_{n+j-1}^2\left(\frac{\int_Af(X_{n+j-1}\mid\theta)dG_{n+j-2}(\theta)}{\int_\Theta f(X_{n+j-1}\mid\theta)dG_{n+j-2}(\theta)} -G_{n+j-2}(A)  \right)^2\\
%&=b_{n+1}\sum_{k\geq n+1}\alpha_{k}^2\left(\frac{\int_Af(X_{k}\mid\theta)dG_{k-1}(\theta)}{\int_\Theta f(X_{k}\mid\theta)dG_{k-1}(\theta)} -G_{k-1}(A)  \right)^2\\
&=b_{n+1}\sum_{k\geq n+1}\frac{Y_{k}}{a_kb_{k}^2},\end{aligned}
\]
where 
\[
Y_k=\left(\frac{\int_Af(X_{k}\mid\theta)dG_{k-1}(\theta)}{\int_\Theta f(X_{k}\mid\theta)dG_{k-1}(\theta)} -G_{k-1}(A)  \right)^2.
\] 
Proceeding as in Lemma \ref{lemma:V}, it can be proved that
\[
E(Y_k\mid\mathcal F_{k-1})=E\left( \left[\frac{\int_Af(X_{k}\mid\theta)dG_{k-1}(\theta)}{\int_\Theta f(X_{k}\mid\theta)dG_{k-1}(\theta)} -G_{k-1}(A)\right]^2\mid\mathcal F_{k-1}  \right)\rightarrow V_A\quad \mbox{$P$-a.s.,}
\]
as $k\rightarrow\infty$. 
Moreover,
$$\sum_{k\geq 1}\frac{E(Y_k^2)}{a_k^2b_{k}^2}<\infty,$$
as $\sum_{k\geq 1} (a_{k}b_{k})^{-2}=\sum_{k\geq 1}\alpha_{k}^4b_{k}^2<\infty
$ and $|Y_k|\leq 1$.
Since, by assumption, 
$
b_{n+1}\sum_{k\geq n+1}(a_kb_k^2)^{-1}=r_n\sum_{k>n}\alpha_k^2\rightarrow 1,
$
then, by Lemma A.1 in \cite{crimaldi2016fluctuation}, 
$U_n\rightarrow V_A$ $P$-a.s. as $n\rightarrow\infty$. 
\\
$\square$

\bigskip
 \noindent {\sc Proof of Theorem \ref{th:central2}}.
  
We first prove that, for every $A\in\mathcal B(\Theta)$, the conditional distribution of $(\sqrt {r_n}\;(G_n(A)-G(A)),{V_{A,n}})$, given $X_{1:n}$, converges to $\Norm(0,{ V_{A}}) \times \delta_{V_A}$, $P$-a.s., on the set $\{\omega : V_A(\omega) >0\}$. 
To show this,  we use Lemma \ref{lemma:V} and compute the joint characteristic function
\[
\begin{aligned}
E(\exp(\mathrm i s_1\sqrt{r_n}\;(G(A)-G_n(A))+\mathrm i s_2{V_{A,n}})\mid 
\cF_{n})&=E(\exp(\mathrm is_1\sqrt{r_n}\;(G(A)-G_n(A)))\mid \cF_{n})\exp(\mathrm is_2{V_{A,n}})\\
&\rightarrow \exp(-s_1^2{ V_{A}}/2)\exp(\mathrm is_2{V_{A}}).
\end{aligned}
\]
Let now $\mathcal D$ be a countable convergence-determining class of bounded continuous functions for the probability measures on $\mathbb R$ and let 
\[
D_n=\sqrt{r_n}\;(G(A)-G_n(A)),\quad 
W_n=\frac{1}{\sqrt{V_{A,n}}}1_{(V_{A,n}>0)},
\quad W=\frac{1}{\sqrt {V_{A}}}1_{(V_A>0)}.
\] 
%Notice that $V_A>0$ if $G(A)\in (0,1)$. 
Then $ W_n(\omega)\rightarrow W(\omega)$ for $P$-almost all $\omega$ such that $V_A(\omega)>0$. 
By Theorem \ref{th:rategeneral}, for every $h\in\mathcal D$,  
$$
E(h(D_n)\mid \cF_{n}) 1_{(V_A>0)} \rightarrow \int h(x) \phi(x \mid 0,1/W^2)dx \quad \mbox{$P$-a.s.},
$$
where $\phi(x\mid \mu, \sigma^2)$ denotes the $\Norm(\mu,\sigma^2)$ density computed at $x$.
Since $W_n$ is  a function of $X_{1:n}$,  then for every $h_1,h_2\in\mathcal D$, 
\[
\begin{aligned}
E(h_1(D_n)h_2(W_n)\mid \cF_{n}) 1_{(V_A>0)} &= E(h_1(D_n)\mid \cF_{n}) h_2(W_n) 1_{(V_A>0)}\\
&\rightarrow 
\int h_1(x) \phi(x \mid 0,1/W^2) dx\, h_2(W) 1_{(V_A>0)}\\
&=\int h_1(x_1) h_2(x_2) d(\Norm(0,1/W^2) \times \delta_W)(x_1, x_2) \, 1_{(V_A>0)}.
\end{aligned}
\]
Since the class $\{h_1 h_2:h_1,h_2\in\mathcal D\}$ is a convergence determining class for the probability measures on ${\mathbb R}^2$, then, for every bounded continuous function $h$,  
\[
\begin{aligned}
E\left(h\left(\sqrt{r_n}\;\frac{G(A)-G_n(A)}{\sqrt { V_{A,n}}} \right)\mid \cF_{n}\right) &=E(h(D_n W_n)\mid \cF_{n}) \\
& \rightarrow \int h(x W) \phi(x \mid 0,1/W^2) dx \\
&=\int h(y) \phi(y \mid 0,1) dy,
\end{aligned}
\]
$P$-a.s. on the set $\{ \omega : V_A(\omega)>0\}$.
$\square$

\bigskip

\noindent {\sc Proof of Theorem \ref{th:ratevett1}}
 
Let $c_1,\dots,c_k$  be arbitrary real numbers. The sequence  
$(\sum_{i=1}^k c_i G_n(A_i))_{n\geq 1}, 
$ is a bounded martingale,  converging  to
$\sum_{i=1}^k c_iG(A_i)$, $P$-a.s.
Following the same steps as in Theorem \ref{th:rategeneral}, with $\sum_{i=1}^k c_i G_n(A_i)$ in the place of $G_n(A)$ and $\sum_{i=1}^k c_iG(A_i)$ in the place of $G(A)$,
 we obtain 
$$ P( \sqrt{r_n}\; (\sum_{i=1}^k c_iG(A_i)-\sum_{i=1}^k c_iG_n(A_i)) \leq t \mid x_{1:n}) \rightarrow \Phi(t ; 0, U), \quad \mbox{for any $t$},$$
where
$U$ is the $P$-a.s. limit of
\[
U_n\equiv r_n \sum_{j\geq n+1} \alpha_j ^2Y_{j},
\]
with
\[
\begin{aligned}
Y_{j}&=\left[\sum_{i=1}^k c_i\left(
P(\theta_{j}\in A_i \mid \cF_{j})-P(\theta_{j}\in A_i\mid \cF_{j-1})	
\right)  \right]^2.
\end{aligned}
\]
Applying Lemma A.1 in \cite{crimaldi2016fluctuation}, as in the proof of Theorem \ref{th:rategeneral}, and noticing that
\[
\begin{aligned}
E(Y_{j}\mid \cF_{j-1})&=\sum_{i,i'} c_i c_{i'} \left[\int_{\mathbb X} P_{G_{j-1}}(A_i \mid x) P_{G_{j-1}}(A_{i'} \mid x)dF_{G_{j-1}}(x)-G_{j-1}(A_i)G_{j-1}(A_{i'})\right]\\
&\rightarrow
\sum_{i,i'} c_ic_{i'}\left[\int_{\mathbb X} P_{G}(A_i \mid x)P_{G}(A_{i'} \mid x)dF_{G}(x)-G(A_i)G(A_{i'})\right]\\
&=\sum_{i,i'=1}^jc_ic_{i'}C_{A_i,A_{i'}}, 
\end{aligned}
\]
we obtain $U=\sum_{i,i'=1}^k c_ic_{i'} C_{A_i,A_{i'}}$. 
Thus, for every $c_1,\dots,c_k$, 
$$ P( \sqrt{r_n} \;(\sum_{i=1}^k c_iG(A_i)-\sum_{i=1}^k c_i G_n(A_i)) \leq t \mid \cF_{n}) \rightarrow \Phi(t  \mid  0,\sum_{i,i'=1}^k c_i c_{i'} C_{A_i,A_{i'}} ), \quad  \mbox{$P$-a.s.,}$$ 
for every $t$. The thesis follows from Cram\'er-Wold theorem.\\
$\square$

\bigskip

 \noindent {\sc Proof of Theorem \ref{th:ratevett2}}.

 Consider the $\cF_n$-measurable spectral decomposition
\[
\frac{ C_n(A_1,\dots,A_k)}{r_n}=Q_n\Lambda_n Q_n^T
\]
where $Q_n$ in a $k\times k$ orthogonal matrix and $\Lambda_n=diag(\lambda_1^{(n)},\dots,\lambda_k^{(n)})$. 
Let $$Y^{(n)}=Q_n^T\left[
\begin{array}{c}
G(A_1)-G_n(A_1)\\
\dots\\
G(A_k)-G_n(A_k)
\end{array}
\right]
$$
and
$$
Z_i^{(n)}=\frac{Y_i^{(n)}}{\sqrt {\lambda_i^{(n)}}}1_{(\lambda_i^{(n)}\neq 0)}+\tilde Z_i 1_{(\lambda_i^{(n)}=0)},
$$
where $\tilde Z_1^{(n)},\dots,\tilde  Z_k^{(n)}$ are i.i.d. random variables, independent of $\mathcal F_n$, and with $\Norm(0,1)$ distribution.  
Then 
$$
\left[
\begin{array}{c}
G(A_1)-G_n(A_1)\\
\dots\\
G(A_k)-G_n(A_k)
\end{array}
\right]=Q_n\Lambda_n^{1/2}Z^{(n)}=\frac{ C_n(A_1,\dots,A_j)^{1/2}}{\sqrt r_n} Z_*^{(n)},
$$
where $ Z_*^{(n)}=Q_nZ^{(n)}\approx \Norm(0, I)$, given $\mathcal F_n$. 

\noindent $\square$

\bigskip

\noindent {\sc Proof of Proposition \ref{th:crebibleRegion} }.

With the same notation as in the proof of Theorem \ref{th:ratevett1}, we can write
\[\begin{aligned}
&\liminf_n P((G(A_1),\dots,G(A_k)\in E_n^{(\epsilon)}\mid\mathcal F_n)\\
&\geq\liminf_n
P\left(Z_*^{(n)T} \frac{ C_n(A_1,\dots,A_k)^{1/2}}
{\sqrt{r_n}}\;\left(\frac{ C_n(A_1,\dots,A_k)+\epsilon I}{r_n}\right)^{-1}\frac{C_n(A_1,\dots,A_k)^{1/2}}{\sqrt{r_n}}Z_*^{(n)}\leq \chi^2_{1-\gamma}\mid\mathcal F_n\right)
\\
%&\geq\liminf_n
%P\left(Z_*^{(n)T}\frac{(C_n(A_1,\dots,A_k)+\epsilon I)^{1/2}}{\sqrt n}\left(\frac{C_n(A_1,\dots,A_k)+\epsilon I}{n}\right)^{-1}\frac{(C_n(A_1,\dots,A_k)+\epsilon I)^{1/2}}{\sqrt n}Z_*^{(n)}\leq \chi^2_{1-\gamma}\mid\mathcal F_n\right)
%\\
&
\geq \liminf_n
P(Z_*^{(n)T}Z_*^{(n)}\leq \chi^2_{1-\gamma}\mid\mathcal F_n)=1-\gamma.
\end{aligned}
\]
\noindent $\square$

\bigskip

{\sc Proof of Proposition \ref{prop:statespace}}. 

For every $A$, $P(Y_2\in A)=E(\tilde H_2(A))=E(\tilde H_1(A))=P(Y_1\in A)$. Moreover, for every $n\geq 2$,
\[
P(Y_{n+1}\in A\mid Y_{1:n-1})
%=E(P(Y_{n+1}\in A\mid \tilde H_{n+1},Y_{1:n})\mid Y_{1:n-1})
=E(\tilde H_{n+1}(A)\mid Y_{1:n-1})
=E(\tilde H_{n}(A)\mid Y_{1:n-1})
%=E(P(Y_{n}\in A\mid \tilde H_{n},Y_{1:n-1})\mid Y_{1:n-1})
=P(Y_{n}\in A\mid Y_{1:n-1}).
\]
Hence $(Y_n)$ is c.i.d. By the properties of c.i.d. sequences, the directing random measure, $\tilde H$, satisfies, $P$-a.s.,
\[
\tilde H(A)=\lim_n P(Y_n\in A\mid Y_{1:n-1})=\lim_n E( P(Y_n\in A\mid Y_{1:n-1},\tilde H_n)\mid Y_{1:n-1})=\lim_n E(\tilde H_n(A)\mid Y_{1:n-1}).\]
\noindent $\square$

\end{document}